\documentclass{article}
\usepackage{amsfonts}
\usepackage{amsmath}
\usepackage{amssymb}
\usepackage{amscd}
\usepackage{amsthm}
\usepackage{hyperref}

\newtheorem{theorem}{Theorem}[section]
\newtheorem{corollary}[theorem]{Corollary}
\newtheorem{lemma}[theorem]{Lemma}

\numberwithin{equation}{section}
\numberwithin{theorem}{section}

\begin{document}

\begin{flushright}
math.SG/0610005
\end{flushright}

\begin{center}
{\LARGE Unitarity in \textquotedblleft quantization commutes with
reduction\textquotedblright} \\[0pt]\bigskip Brian C. Hall\footnote{E-mail:
bhall@nd.edu} \footnote{Supported in part by NSF Grant DMS-02000649.} and
William D. Kirwin\footnote{E-mail: wkirwin@nd.edu} \\[0pt]\bigskip\emph{Dept.
of Mathematics, University of Notre Dame, Notre Dame, IN 46556-4618}
\end{center}

\begin{abstract}
Let $M$\ be a compact K\"{a}hler manifold equipped with a Hamiltonian action
of a compact Lie group $G$. In this paper, we study the geometric quantization
of the symplectic quotient $M/\!\!/G$. Guillemin and Sternberg [Invent. Math.
\textbf{67} (1982), 515--538] have shown, under suitable regularity
assumptions, that there is a natural invertible map between the quantum
Hilbert space over $M/\!\!/G$ and the $G$-invariant subspace of the quantum
Hilbert space over $M.$

Reproducing other recent results in the literature, we prove that in general
the natural map of Guillemin and Sternberg is not unitary, \textit{even to
leading order in Planck's constant}. We then modify the quantization procedure
by the \textquotedblleft metaplectic correction\textquotedblright\ and show
that in this setting there is still a natural invertible map between the
Hilbert space over $M/\!\!/G$ and the $G$-invariant subspace of the Hilbert
space over $M.$ We then prove that this modified Guillemin--Sternberg map is
asymptotically unitary to leading order in Planck's constant. The analysis
also shows a good asymptotic relationship between Toeplitz operators on $M$
and on $M/\!\!/G.$

\end{abstract}


\section{Introduction}

Let $M$ be an integral compact K\"{a}hler manifold with symplectic form
$\omega$. Following the program of geometric quantization, suppose we are
given a Hermitian holomorphic line bundle $\ell$ with connection over $M,$
chosen in such a way that the curvature of $\ell$ is equal to $-i\omega.$ We
then consider $\ell^{\otimes k},$ the $k$th tensor power of $\ell,$ where in
this setting we interpret $k$ as the reciprocal of Planck's constant $\hbar. $
The Hilbert space of quantum states associated to $M$ is then (for a fixed
value of $\hbar=1/k$) the space of holomorphic sections of $\ell^{\otimes k}.
$

Now suppose that we are given a Hamiltonian action of a connected compact Lie
group $G$ on $M.$ Then we can construct the symplectic quotient $M/\!\!/G,$
which is another compact K\"{a}hler manifold (under suitable regularity
assumptions on the action of $G$). The line bundle $\ell$ naturally descends
to a bundle $\hat{\ell}$ over $M/\!\!/G,$ and the quantum Hilbert space
associated to $M/\!\!/G$ is then the space of holomorphic sections of
$\hat{\ell}^{\otimes k}.$ This space of sections is what we obtain by
performing reduction by $G$ \textit{before} quantization. Alternatively, we
may \textit{first} perform quantization and \textit{then} perform reduction at
the quantum level, which amounts to restricting to the space of $G$-invariant
sections of $\ell^{\otimes k}.$

A classic result of Guillemin and Sternberg is that there is a natural
invertible linear map between the \textquotedblleft first reduce and then
quantize\textquotedblright\ space (the space of all holomorphic sections of
the bundle $\hat{\ell}^{\otimes k}$ over $M/\!\!/G$) and the \textquotedblleft
first quantize and then reduce\textquotedblright\ space (the space of
$G$-invariant holomorphic sections of the bundle $\ell^{\otimes k}$ over $M$).
This result is sometimes described as saying that quantization commutes with
reduction. However, from the point of view of quantum mechanics, it is not
just the vector space structure of the quantum Hilbert space that is
important, but also the inner product.

It is natural, then, to investigate the extent to which the
Guillemin--Sternberg map is unitary. Some work in this direction has already
been done, and in general, it has been found that the obstruction to unitarity
is a certain function on $M/\!\!/G$ (sometimes referred to as the effective
potential). This obstruction, in some form, was identified independently in
the works of \cite{Flude}, \cite{Charles}, and \cite{Paoletti}, and
\cite{Ma-Zhang06}. See also \cite{Gotay}, \cite{HallGQ}, and
\cite{Huebschmann}. We discuss the relation of our work to these in Section
\ref{prior.sec}.

Our first main result is a new proof (similar to that of Charles in
\cite{Charles} for the torus case) that the Guillemin--Sternberg map is not
unitary in general, and indeed that the map does not become asymptotically
unitary as $k\rightarrow\infty$ (i.e., as $\hbar\rightarrow0$). We show that
the obstruction to asymptotic unitarity of the map is the volume of the
$G$-orbits inside the zero-set of the moment map in $M.$ If these $G$-orbits
do not all have the same volume, then the Guillemin--Sternberg map will not be
asymptotically unitary. This failure of asymptotic unitarity is troubling from
a physical point of view. After all, although one expects different
quantization procedures (e.g., performing quantization and reduction in
different orders) to give different results, one generally regards these
differences as \textquotedblleft quantum corrections\textquotedblright\ that
should disappear as $\hbar$ tends to zero.

The main contribution of this paper is a solution to this unsatisfactory
situation: we include the so-called metaplectic correction, which involves
tensoring the original line bundle with the square root of the canonical
bundle (assuming that such a square root exists). We first show that one can
define a natural map of Guillemin--Sternberg type in the presence of the
metaplectic correction, and this natural map is invertible for all
sufficiently large values of the tensor power $k.$ We then show that this
modified Guillemin--Sternberg map, unlike the original one, is asymptotically
unitary in the limit as $k$ tends to infinity.

In the rest of this introduction, we describe our results in greater detail
and compare them to previous results of other authors.

\subsection{Main results}

Let $M^{2n}$ be a compact K\"{a}hler manifold with symplectic form $\omega$,
complex structure $J,$ and Riemannian metric $B=\omega(\cdot,J\cdot).$ Assume
that $M$ is quantizable, that is, that the class $[\omega/2\pi]$ is integral.
Choose a Hermitian line bundle $\ell$ over $M$ with compatible connection
$\nabla$ in such a way that the curvature of $\nabla$ is $-i\omega.$ (Such a
bundle exists because $M$ is quantizable.) We call $\ell$ a prequantum bundle
for $M.$ We denote the Hermitian form on $\ell,$ which we take to be linear in
the second factor, by $(\cdot,\cdot)$ and we denote the pointwise magnitude of
a section $s$ of $\ell$ by $\left\vert s\right\vert ^{2}(x)=(s,s)(x).$ The
connection and Hermitian form on $\ell$ induce a connection and Hermitian form
on $\ell^{\otimes k}$ which we denote by the same symbols.

The complex structure on $M$ gives $\ell^{\otimes k}$ the structure of a
holomorphic line bundle in which the holomorphic sections of are those that
are covariantly constant in the $(0,1)$-directions. We let $\mathcal{H}%
(M;\ell^{\otimes k})$ denote the (finite-dimensional) space of holomorphic
sections of $\ell^{\otimes k}.$ The symplectic form on $M$ induces a volume
form
\begin{equation}
\varepsilon_{\omega}:=\frac{\omega^{\wedge n}}{n!},\label{liouville}%
\end{equation}
which is known as the Liouville volume form. We make $\mathcal{H}%
(M;\ell^{\otimes k})$ into a Hilbert space by considering the inner product
given by%
\begin{equation}
\langle s_{1},s_{2}\rangle:=\left(  k/2\pi\right)  ^{n/2}\int_{M}(s_{1}%
,s_{2})\,\varepsilon_{\omega}.\label{innerproduct}%
\end{equation}

In geometric quantization, we interpret the tensor power $k$ as the reciprocal
of Planck's constant $\hbar.$ Thus the study of holomorphic sections of
$\ell^{\otimes k}$ in the limit $k\rightarrow\infty,$ familiar from algebraic
geometry, is in this setting interpreted as the \textquotedblleft
semiclassical limit\textquotedblright\ (i.e., the limit $\hbar\rightarrow0$).
The semiclassical limit in geometric quantization has garnered much recent
interest. See for example the work of Borthwick and Uribe \cite{UribeGQ} and
\cite{Uribe2} on a symplectic version of the Kodaira embedding theorem and the
work of Borthwick, Paul, and Uribe \cite{Uribe3} on applications to relative
Poincar\'{e} series.

Suppose now that we are given a Hamiltonian action of a connected compact Lie
group $G$ of dimension $d,$ together with an equivariant moment map
$\mathrm{\Phi}$ for this action. We assume (as in \cite{Guillemin-Sternberg})
that $0$ is in the image of $\Phi,$ that $0$ is a regular value of $\Phi,$ and
that $G$ acts freely on $\Phi^{-1}(0).$ The symplectic quotient $M/\!\!/G $ is
then defined to be
\[
M/\!\!/G=\Phi^{-1}(0)/G.
\]
The quotient $M/\!\!/G$ acquires from $M$ the structure of a quantizable
K\"{a}hler manifold of dimension $2(n-d)$; in particular the symplectic form
$\omega$ on $M$ descends to a symplectic form $\widehat{\omega}\in\Omega
^{2}(M/\!\!/G)$. \ We assume that the action of $G$ lifts to $\ell.$ Then
$\ell$ descends to a holomorphic Hermitian line bundle $\hat{\ell}$ with
connection over $M/\!\!/G,$ with curvature $-i\widehat{\omega}.$ Let
$\mathcal{H}(M/\!\!/G;\hat{\ell}^{\otimes k})$ denote the space of holomorphic
sections of $\hat{\ell}^{\otimes k}$ and use on this space the inner product
given by%
\begin{equation}
\langle s_{1},s_{2}\rangle:=\left(  k/2\pi\right)  ^{(n-d)/2}\int
_{M/\!\!/G}(s_{1},s_{2})~\varepsilon_{\widehat{\omega}}.\label{innerproduct2}%
\end{equation}

We are interested in the relationship between two different Hilbert spaces:
first, the space of $G$-invariant holomorphic sections of $\ell^{\otimes k}$
over $M,$ with the inner product (\ref{innerproduct}); and second, the space
of all holomorphic sections of $\hat{\ell}^{\otimes k}$ over $M/\!\!/G,$ with
the inner product (\ref{innerproduct2}). The first Hilbert space, denoted
$\mathcal{H}(M;\ell^{\otimes k})^{G},$ is the one obtained by first quantizing
and then reducing by $G$; the second one is obtained by first reducing by $G$
and then quantizing. According to Guillemin and Sternberg
\cite{Guillemin-Sternberg}, there is (for each $k$) a natural one-to-one and
onto linear map $A_{k}$ from $\mathcal{H}(M;\ell^{\otimes k})^{G}$ onto
$\mathcal{H}(M/\!\!/G;\hat{\ell}^{\otimes k}).$ In particular, these two
spaces have the same dimension.

The map $A_{k}$ is defined in the only reasonable way: one takes a
$G$-invariant holomorphic section $s$ of $\ell^{\otimes k},$ restricts it to
$\Phi^{-1}(0),$ and then lets $s$ descend from $\Phi^{-1}(0)$ to $\Phi
^{-1}(0)/G.$ From the way that the complex structure on $M/\!\!/G$ is defined,
it is easy to see that $A_{k}$ maps holomorphic sections of $\ell^{\otimes k}$ to
holomorphic sections of $\hat{\ell}^{\otimes k}.$ The hard part is to show
that this map is invertible, i.e., that a section of $\hat{\ell}^{\otimes k},$
after being lifted to $\Phi^{-1}(0),$ can be extended holomorphically to a
($G$-invariant) section over all of $M.$ In proving this, Guillemin and
Sternberg make use of a holomorphic action of the \textquotedblleft
complexified\textquotedblright\ group $G_{\mathbb{C}}$ on $M.$

We wish to determine the extent to which the map $A_{k}$ is unitary. We show
that (for each $k$) there exists a function $I_{k}$ on $M/\!\!/G$ with the
property that for each $G$-invariant holomorphic section $s$ of $\ell^{\otimes
k}$ we have%
\[
\left\Vert s\right\Vert ^{2}=\left(  k/2\pi\right)  ^{n/2}\int_{M}\left\vert
s\right\vert ^{2}\varepsilon_{\omega}=\left(  k/2\pi\right)  ^{(n-d)/2}%
\int_{M/\!\!/G}\left\vert A_{k}s\right\vert ^{2}I_{k\,}\varepsilon
_{\widehat{\omega}}.
\]
Our first main result is the following.

\begin{theorem}
\label{lim1.thm}The functions $I_{k}$ satisfy%
\[
\lim_{k\rightarrow\infty}I_{k}([x_{0}])=2^{-d/2}\mathrm{\operatorname*{vol}%
}(G\cdot x_{0})
\]
for all $x_{0}\in\Phi^{-1}(0),$ and the limit is uniform. Here
\textrm{$\operatorname*{vol}$}$(G\cdot x_{0})$ denotes the volume of the
$G$-orbit of $x_{0} $ with respect to the Riemannian structure inherited from
$M.$
\end{theorem}

If all of the $G$-orbits in the zero-set of the moment map have the same
volume, then Theorem \ref{thm:density asymptotics} implies that the natural
map $A_{k}$ is asymptotically a constant multiple of a unitary map in the
limit as $\hbar$ tends to zero. If, as is usually the case, the $G$-orbits in
the zero-set do not all have the same volume, then using peaked sections, we
can show that $A_{k}$ is not asymptotically a constant multiple of a unitary
map, in a sense described in Theorem \ref{thm:non-unitarity}.

Now, it is probably unreasonable to expect the natural map between the
\textquotedblleft first reduce and then quantize\textquotedblright\ space and
the \textquotedblleft first quantize and then reduce\textquotedblright\ space
to be unitary. Different ways of performing quantization (e.g., before
reduction and after reduction) in general give inequivalent results,
especially if equivalence is measured by the existence of a geometrically
natural \textit{unitary} (as opposed to merely invertible) map. On the other
hand, one expects differences between the procedures to be \textquotedblleft
quantum corrections\textquotedblright\ that vanish for small $\hbar.$ The
situation reflected in Theorem \ref{lim1.thm}, then, is unsatisfactory; the
geometrically natural map between $\mathcal{H}(M;\ell^{\otimes k})^{G}$ and
$\mathcal{H}(M/\!\!/G;\hat{\ell}^{\otimes k})$ is not unitary \textit{even to
leading order in }$\hbar.$

To remedy this situation, we introduce the \textquotedblleft
metaplectic\textquotedblright\ or \textquotedblleft
half-form\textquotedblright\ correction. There are various results which seem
to indicate the geometric quantization works better with half-forms (indeed,
our main results can be seen as further justification of this claim). For
example, it is only with half-forms that the quantization of the simple
harmonic oscillator has the correct zero-point energy (see \cite[Chap. 10]{Woodhouse}
for details and further examples). For this reason, we might expect (and it is
indeed the case) that the inclusion of half-forms remedies the unsatisfactory
situation described by Theorem \ref{lim1.thm}.

The metaplectic correction consists of tensoring $\ell^{\otimes k}$ with a
square root of the canonical bundle $K$ of $M$, assuming such a square root
exists (equivalently, the vanishing of the second Stiefel--Whitney class), and
similarly tensoring $\hat{\ell}^{\otimes k}$ with a square root of
the canonical bundle $\widehat{K}$ of $M/\!\!/G.$ Assuming that the $G$-action
lifts to $\ell^{\otimes k}\otimes\sqrt{K}$, our second main result is an analog of
\cite[Thm 5.2]{Guillemin-Sternberg} in the presence of half-forms.

\begin{theorem}
\label{bk.thm}For each $k,$ there is a natural linear map $B_{k}$ between the
space of $G$-invariant holomorphic sections of $\ell^{\otimes k}\otimes
\sqrt{K}$ over $M$ and the space of all holomorphic sections of $\hat{\ell
}^{\otimes k}\otimes\sqrt{\widehat{K}}$ over $M/\!\!/G.$ Furthermore, this map
is invertible for all sufficiently large $k.$
\end{theorem}

We then show that there exists a function $J_{k}$ on $M/\!\!/G$ such that for
each $G$-invariant holomorphic section $r$ of $\ell^{\otimes k}\otimes\sqrt
{K}$ we have%
\[
\left\Vert r\right\Vert ^{2}:=\left(  k/2\pi\right)  ^{n/2}\int_{M}\left\vert
r\right\vert ^{2}\varepsilon_{\omega}=\left(  k/2\pi\right)  ^{(n-d)/2}%
\int_{M/\!\!/G}\left\vert B_{k}r\right\vert ^{2}J_{k}\varepsilon
_{\widehat{\omega}}.
\]
Our last main result is the following, which implies that the maps $B_{k},$
unlike the maps $A_{k},$ are asymptotically unitary in the limit
$k\rightarrow\infty$ (i.e., $\hbar\rightarrow0$).

\begin{theorem}
\label{lim2.thm}The functions $J_{k}$ satisfy%
\[
\lim_{k\rightarrow\infty}J_{k}([x_{0}])=1
\]
for all $x_{0}\in\Phi^{-1}(0),$ and the limit is uniform.
\end{theorem}

The origin of the unwanted volume factor in Theorem \ref{lim1.thm} is the
relationship between the volume measure on $\Phi^{-1}(0)$ and the volume
measure on $\Phi^{-1}(0)/G$ (Lemma \ref{lemma:zero-set integral decomp}). The
metaplectic correction introduces a compensating volume factor in the
pointwise magnitude of a section over $\Phi^{-1}(0)$ and the pointwise
magnitude of the corresponding section over $\Phi^{-1}(0)/G$ (Theorem
\ref{thm:BKS comparison}).

The proofs of our main results make use of a holomorphic action of the
\textquotedblleft complexified\textquotedblright\ group $G_{\mathbb{C}}$ on
$M$, obtained in \cite{Guillemin-Sternberg} by analytic continuation of the
action of $G.$ We let $M_{s}$ denote the \textbf{stable set,} that is, the set
of points in $M$ that can be moved into the zero-set of the moment map by the
action of $G_{\mathbb{C}}.$ The stable set is an open set of full measure in
$M$ and the symplectic quotient $M/\!\!/G=\Phi^{-1}(0)/G$ can be identified
naturally with $M_{s}/G_{\mathbb{C}}.$ We show how the magnitude of a
$G$-invariant section (of $\ell^{\otimes k}$ or $\ell^{\otimes k}\otimes
\sqrt{K}$) varies in a predictable way as one moves off of $\Phi^{-1}(0)$
along a $G_{\mathbb{C}}$-orbit; this fact, which is crucial to our analysis,
has been noticed before by Guillemin and Sternberg \cite{Guillemin-Sternberg}
and by Donaldson \cite{Donaldson}. The densities $I_{k} $ and $J_{k}$ are
obtained by integrating certain quantities related to the moment map over each
$G_{\mathbb{C}}$-orbit. The asymptotic behavior of these densities is then
determined by Laplace's method (sometimes referred to as the stationary phase
approximation\ or the method of steepest descent).

In proving the invertibility of our Guillemin--Sternberg-type map $B_{k},$ we
first show that a $G$-invariant section of $\ell^{\otimes k}\otimes\sqrt{K}$
defined over $\Phi^{-1}(0)$ can be extended to a holomorphic, $G_{\mathbb{C}}%
$-invariant section over the stable set $M_{s}.$ We then show that \textit{for
all sufficiently large }$k,$ these sections have a removable singularity over
the complement of the stable set and thus extend to holomorphic sections over
all of $M.$ Here, the presence of half-forms makes the situation slightly more
complicated than in \cite{Guillemin-Sternberg}, where the natural map $A_{k}$
is shown to be invertible for all $k.$

A trivial modification of the arguments described above yields a formula that
relates a Toeplitz operator with a $G$-invariant symbol upstairs, restricted
to the space of $G$-invariant sections, to a certain Toeplitz operator
downstairs. Suppose $f$ is a smooth, $G$-invariant function on $M.$ Then in
the case with half-forms, the result is that the Toeplitz operator with symbol
$f$ upstairs, restricted to the $G$-invariant subspace, is asymptotically
equivalent to the Toeplitz operator with symbol $\hat{f}$ downstairs, where
$\hat{f}$ is obtained by restricting $f$ to $\Phi^{-1}(0)$ and letting it
descend to $\Phi^{-1}(0)/G.$

\bigskip

For the reader's convenience, we collect here the various assumptions made in
this paper. We begin with a compact K\"{a}hler manifold $M$ of real dimension
$2n,$ with symplectic form $\omega.$ We assume that $\left[ \omega/2\pi\right]
$ is an integral cohomology class. We assume that $M$ is equipped with a
holomorphic and Hamiltonian action of a compact $d$ -dimensional Lie group $G$
with equivariant moment map $\Phi:M\rightarrow\mathfrak{g}^{\ast}.$ We assume
that $0$ is in the image of the moment map, that $0$ is a regular value of the
moment map, and that $G$ acts freely on the zero-set $\Phi^{-1}(0)$. Since the
symplectic form is integral, there exists over $M$ at least one complex
Hermitian line bundle with compatible connection with curvature $-i\omega$; we
fix some choice and denote it by $\ell$. In a canonical way, the moment map
defines an infinitesimal holomorphic action of $G$ on $\ell$. We assume that
this infinitesimal action can be exponentiated. These assumptions are the same
as in \cite{Guillemin-Sternberg} and are all the assumptions relevant to
Theorem \ref{lim1.thm}.

We need additional assumptions for Theorems \ref{bk.thm} and \ref{lim2.thm}.
Let $K=\bigwedge^{n}(T^{\ast}M)^{1,0}$ denote the canonical bundle of $M$. We
assume that $K$ admits a square root $\sqrt{K}\rightarrow M$ (this is
equivalent to the vanishing of the second Stiefel-Whitney class of $M$), and
finally we assume that the $G$-action on $K$ (which is induced from that on
$M$) can be lifted to an action on $\sqrt{K}.$

\subsection{Prior results\label{prior.sec}}

To the best of our knowledge, the first work that systematically address the
question of unitarity in the context of quantization and reduction of
K\"{a}hler manifolds is the Ph.D. thesis of Flude \cite{Flude}. Flude gives a
formal computation of the leading-order asymptotics of the density $I_{k}$ as
$k$ tends to infinity, without the metaplectic correction. Flude's main
computation is essentially the same as ours (in the case without the
metaplectic correction): an application of Laplace's method. However, Flude's
result is not rigorous, because he considers the magnitude of an invariant
section only in a small neighborhood of the zero-set. This should be
compared to our Theorem \ref{thm:s along a curve}, which gives an expression
for the magnitude of an invariant section that is valid everywhere. (Even
after Theorem \ref{thm:s along a curve} is obtained, some control is required
over the blow-up of certain Jacobians near the unstable set.) Nevertheless,
Flude does identify the volume of the $G$-orbits in $\Phi^{-1}(0)$ as the
obstruction to the asymptotic unitarity of the Guillemin--Sternberg map.

Next, there is the work of L. Charles \cite{Charles}, who considers the case
in which the group $G$ is a torus and without half-forms. (Charles does
consider half-forms in a related context in \cite{Charles2}.) In the torus
case, Proposition 4.17 of \cite{Charles}, with $f$ identically equal to $1$ in
a neighborhood of the zero-set, is essentially the same as our Theorem
\ref{thm:s norm decomp}. Proposition 4.18 of \cite{Charles} is then
essentially the same as the first part of our Theorem
\ref{thm:density asymptotics}. We do, however, go slightly beyond the approach
of Charles, in that we give an \textit{exact} (not just asymptotic to all
orders) expression for the norm of an invariant section upstairs as an
integral over the downstairs manifold. This requires some control over the
blow-up of certain Jacobians as one approaches the unstable set.

Then there is the work of Paoletti \cite{Paoletti}. He takes a different
approach to measuring unitarity, looking at the behavior of orthonormal bases.
Nevertheless, his result agrees with ours (in the case without half-forms) in
that he identifies the volume of the $G$-orbits in the zero-set of the moment
map as an obstruction to unitarity of the Guillemin--Sternberg map.

Finally, there is the the work of Ma and Zhang in \cite{Ma-Zhang06}. In this
paper, the authors compute an asymptotic expansion of the equivariant Bergman
kernel, and our Theorem \ref{lim1.thm} is a special case of their Theorem 0.1.
In contrast to our work and the works of Flude and Charles (all of which are
based on Laplace's method), those of Paoletti and Ma--Zhang are based on the
microlocal analysis developed by L. Boutet de Monvel and Guillemin in
\cite{BdMG}.

To the best of our knowledge, ours is the first paper that constructs a
Guillemin--Sternberg-type map in the presence of the metaplectic correction
and therefore the first to show in a general setting that the metaplectic
correction improves the situation in regard to unitarity. (Such an improvement
had already been suggested by the special examples considered in \cite{DHym}
and \cite{HallGQ}. See the last section for a discussion of these examples.)

We close this section by mentioning the work of Gotay \cite{Gotay}, who
considers the relationship between quantization and reduction in the context
of a cotangent bundle equipped with the vertical polarization (a \emph{real}
polarization). Gotay considers the cotangent bundle of an $n$-manifold $Q$ and
assumes that a compact group $G$ of dimension $d$ acts freely on $Q.$ The
induced action of $G$ on $T^{\ast}Q$ is then free and Hamiltonian. Gotay
includes half-forms in the quantization (as one must in the case of the
vertical polarization) and obtains a \emph{unitary} map between the
quantize-then-reduce space and the reduce-then-quantize space.

Gotay's work does not overlap with ours, because we consider only complex
polarizations. Nevertheless, it is worth noting that Gotay obtains
\textit{exact} unitarity, whereas we obtain unitarity only asymptotically as
Planck's constant tends to zero, even with half-forms. The reason that the
cotangent bundle case is nicer than the K\"{a}hler case seems to be that there
is less differential geometry involved. There is no additional structure on
the manifold $Q$ (metric or measure or complex structure) that enters. Gotay's
unitarity result comes down to the relationship between the integral over $Q$
of a $G$-invariant $n$-form $\alpha$ and the integral over $Q/G$ of an
$(n-d)$-form $\hat{\alpha}$ obtained by contracting $\alpha$ with the
generators of the $G$-action and then letting the result descend to $Q/G.$
(Here $\alpha$ is the square of a half-form.) In the K\"{a}hler case, we have
a mapping between $(n,0)$-forms on $M$ and $(n-d,0)$-forms on $M/\!\!/G$ defined
in a manner very similar to that in \cite{Gotay}. However, it does not make
sense to integrate an $(n,0)$-form over $M,$ because $M$ has real dimension
$2n.$ So the norm of a half-form cannot be computed by simply squaring and
integrating. Rather, one uses a more complicated procedure involving both the
complex structure on $M$ and the Liouville volume measure. The relationship
between the norm upstairs and the norm downstairs is correspondingly more
complicated and involves the geometric structures that we have on $M$ but not
on $Q.$

Another way of thinking about the K\"{a}hler case is to observe that the way
one lifts a section downstairs to a section upstairs is by lifting to a
$G_{\mathbb{C}}$-invariant (not just $G$-invariant) section upstairs. However,
the action of $G_{\mathbb{C}}$ preserves only the complex structure on $M$ and
not the symplectic structure on $M$ or the Hermitian structure on the relevant
line bundles. Thus the relationship between the upstairs norm and the
downstairs norm involves the way the volume measure and the magnitude of an
invariant section vary over each $G_{\mathbb{C}}$-orbit. In the cotangent
bundle case, by contrast, only a $G$-action is involved and this preserves all
the relevant structure.

\section{\label{sec:preliminaries}Preliminaries}

We begin this section by recalling the method of geometric quantization (at
the moment, without half-forms) as it applies to K\"{a}hler manifolds. We then
recall the notion of the Marsden--Weinstein or symplectic quotient and explain
the special form this construction takes in the setting of K\"{a}hler
manifolds. Finally, we describe the natural invertible map, due to Guillemin
and Sternberg, between the \textquotedblleft quantize then
reduce\textquotedblright\ space and the \textquotedblleft reduce then
quantize\textquotedblright\ space.\

\subsection{K\"{a}hler quantization}

Let $(M,\omega,J,B=\omega(\cdot,J\cdot))$ be a K\"{a}hler manifold with
symplectic form $\omega$, complex structure $J$, and Riemannian metric
$B=\omega(\cdot,J\cdot)$. We assume that $M$ is connected, compact, and
integral (i.e., that $[\omega/2\pi]$ is an integral cohomology class). We fix
once and for all a Hermitian line bundle $\ell$ with compatible connection
$\nabla,$ chosen in such a way that the curvature of $\nabla$ is equal to
$-i\omega.$ The connection on $\ell$ induces a connection on the $k$th tensor
power $\ell^{\otimes k},$ also denoted $\nabla.$ For any $k,$ $\ell^{\otimes
k}$ may be given the structure of a holomorphic line bundle in such a way that
the holomorphic sections are precisely those that are covariantly constant in
the $(0,1)$ (or $\bar{z}$) directions.

We interpret the tensor power $k$ as the reciprocal of Planck's constant
$\hbar.$ For each fixed value of $k$ (or $\hbar$), the quantum Hilbert space
is the space of holomorphic sections of $\ell^{\otimes k},$ denoted
$\mathcal{H}(M;\ell^{\otimes k}).$ If $\varepsilon_{\omega}$ denotes the
Liouville volume form (\ref{liouville}), then we use the following natural
inner product on $\mathcal{H}(M;\ell^{\otimes k}),$%
\[
\langle s_{1},s_{2}\rangle:=\left(  k/2\pi\right)  ^{n/2}\int_{M}(s_{1}%
,s_{2})\,\varepsilon_{\omega},
\]
as in (\ref{innerproduct}), where $(s_{1},s_{2})$ is the (pointwise) Hermitian
structure on $\ell^{\otimes k}.$

\subsection{K\"{a}hler reduction\label{subsec:Kahler reduction}}

We now assume that we are given a smooth action of a connected compact Lie
group $G$ of dimension $d$ on $M.$ We assume that the action preserves all of
the structure (symplectic, complex, Riemannian) of $M.$ Let $\mathfrak{g}$
denote the Lie algebra of $G$ and for each $\xi\in\mathfrak{g},$ let $X^{\xi}
$ denote the vector field describing the infinitesimal action of $\xi$ on $M.
$ That is, $X^{\xi}(x)=\frac{d}{dt}\left.  e^{t\xi}\cdot x\right\vert _{t=0}.$

We assume that for each $\xi\in\mathfrak{g}$ there exists a smooth function
$\phi_{\xi}$ on $M$ such that $X^{\xi}$ is the Hamiltonian vector field
associated to $\phi_{\xi};$ i.e., such that
\[
\mathfrak{i}(X^{\xi})\omega=d\phi_{\xi}.
\]
(This is automatically the case if $G$ is semisimple or if $M$ is simply
connected.) For each $\xi,$ the function $\phi_{\xi}$ is unique up to a
constant. Since $M$ is compact, it is always possible to choose the constants
in such a way that the map $\xi\rightarrow\phi_{\xi}$ is linear and that
$\{\phi_{\xi},\phi_{\eta}\}=-\phi_{\lbrack\xi,\eta]}, $ and we fix one
particular choice of constants with these two properties. (For any $G,$ one
way to choose the constants with these two properties is to require each
$\phi_{\xi}$ to have integral zero over $M.$ If $G$ is semisimple, then these
two properties uniquely determine the choice of the constants.)

We may put together the functions $\phi_{\xi}$ into a \textquotedblleft moment
map\textquotedblright\ $\Phi:M\rightarrow\mathfrak{g}^{\ast}$ given by
\[
\Phi(x)(\xi)=\phi_{\xi}(m)
\]
for each $\xi\in\mathfrak{g}$ and $m\in M.$ The condition $\{\phi_{\xi}%
,\phi_{\eta}\}=-\phi_{\lbrack\xi,\eta]}$ ensures that the moment map is
equivariant with respect to the action of $G$ on $M$ and the coadjoint action
of $G$ on $\mathfrak{g}^{\ast}.$ We assume that $0$ is in the image of the
moment map and that it is a regular value, so that the zero-set $\Phi^{-1}(0)$
is a submanifold of $M$. We assume, moreover, that $G$ acts freely on
$\Phi^{-1}(0)$ so that the quotient is a manifold. The quotient%
\[
M/\!\!/G:=\Phi^{-1}(0)/G
\]
is then called the symplectic or Marsden--Weinstein \cite{Marsden-Weinstein}
quotient of $M$ by $G$.

The quotient $M/\!\!/G$ inherits a symplectic structure from $M$: there is a
unique symplectic form $\widehat{\omega}\in\Omega^{2}(M/\!\!/G)$ such that
$i^{\ast}\omega=\pi^{\ast}\widehat{\omega}$, where $i:\Phi^{-1}(0)\rightarrow
M$ denotes the inclusion map and $\pi:\Phi^{-1}(0)\rightarrow\Phi^{-1}(0)/G$
is the quotient map. Furthermore, $[\widehat{\omega}/2\pi]$ is an
integral cohomology class on $M/\!\!/G.$

So far, we have described how the symplectic structure on $M$ induces a
symplectic structure on $M/\!\!/G.$ We now show how the complex structure on
$M$ induces a complex structure on $M/\!\!/G.$ The descent of the complex
structure can be understood either \textquotedblleft
infinitesimally\textquotedblright\ (by describing how the distribution of
$(1,0)$-vectors descends) or \textquotedblleft globally\textquotedblright\ (by
realizing $M/\!\!/G$ as the quotient of the \textquotedblleft stable
set\textquotedblright\ in $M$ by the complexification of $G$).

We begin with the infinitesimal approach. Let $T^{1,0}M$ be the K\"{a}hler
polarization on $M$; i.e.,$T^{1,0}M$ is the $n$-dimensional (complex)
distribution consisting of type-$(1,0)$ vector fields on $M$. For future use,
we note that the projection $\pi_{+}:T^{\mathbb{C}}M\rightarrow T^{1,0}M $ is
given by%
\[
\pi_{+}X=\frac{1}{2}(1-iJ)X.
\]

In \cite{Guillemin-Sternberg}, Guillemin and Sternberg show that the
$G$-orbits in the zero-set are totally real submanifolds; i.e.,for all
$x_{0}\in\Phi^{-1}(0),$%
\[
T_{x_{0}}^{1,0}M\cap\{X_{x_{0}}^{\xi}:\xi\in\mathfrak{g}\}=\{0\}.
\]
Moreover, they show that $T^{\mathbb{C}}(\Phi^{-1}(0))\cap T^{1,0}M$ is a $G
$-invariant complex distribution of complex rank $n-d$, and that $\pi_{\ast
}(T^{\mathbb{C}}(\Phi^{-1}(0))\cap T^{1,0}M)$ is a well-defined integrable
complex distribution of rank $n-d$ on $M/\!\!/G$. In fact, the complex
structure induced by this distribution defines a K\"{a}hler structure on
$M/\!\!/G$, and we henceforth make the identification
\begin{equation}
T^{1,0}\left(  M/\!\!/G\right)  =\pi_{\ast}(T^{\mathbb{C}}(\Phi^{-1}(0))\cap
T^{1,0}M).\label{infinitesimal complex}%
\end{equation}

We now turn to the global approach, which is more intuitive and much more
useful for the computations we will carry out in this paper. We have assumed
that the action of $G$ is holomorphic. Using this, it can be shown
\cite{Guillemin-Sternberg} that the action of $G$ can be analytically
continued to a holomorphic action of the \textquotedblleft
complexified\textquotedblright\ group $G_{\mathbb{C}}$ on $M;$ the
infinitesimal action of this continuation is defined by%
\begin{equation}
X^{i\xi}:=JX^{\xi},~\xi\in\mathfrak{g.}\label{eqn:inf GC action}%
\end{equation}

Here $G_{\mathbb{C}}$ is a connected complex Lie group containing $G$ as a
maximal compact subgroup, and the Lie algebra $\mathfrak{g}_{\mathbb{C}}$ of
$G_{\mathbb{C}}$ is the complexification of $\mathfrak{g}.$ The Cartan
decomposition is a diffeomorphism $G_{\mathbb{C}}\simeq\exp(i\mathfrak{g})G.$
Moreover, the set $\exp(i\mathfrak{g})$ is diffeomorphic to the vector space
$\mathfrak{g}$, the diffeomorphism being the exponential map, and so as smooth
manifolds we have $G_{\mathbb{C}}=\mathfrak{g}\times G$. See \cite[Sec.
6.3]{Knapp} for details.

We let the \textbf{stable set} $M_{s}$ denote the saturation of the zero-set
$\Phi^{-1}(0)$ by the action of $G_{\mathbb{C}}$:
\begin{equation}
M_{s}:=G_{\mathbb{C}}\cdot\Phi^{-1}(0).\label{eqn:stable def}%
\end{equation}
That is, $M_{s}$ is the set of points in $M$ that can be moved into $\Phi
^{-1}(0)$ by the action of $G_{\mathbb{C}}.$ (See
(\ref{eqn:stable set character}) below for another characterization of the
stable set.) The following properties of $M_{s}$ follow from results in
\cite{Guillemin-Sternberg}: (1) $M_{s}$ is an open set of full measure in $M$;
(2) $G_{\mathbb{C}}$ acts freely on $M_{s}$; and (3) each $G_{\mathbb{C}}%
$-orbit in $M_{s}$ intersects $\Phi^{-1}(0)$ in precisely one $G$-orbit.

We will show shortly (Theorem \ref{thm:gc bundle}) that $M_{s}$ is in fact a
principal $G_{\mathbb{C}}$-bundle over $\Phi^{-1}(0)/G.$ This implies
(Corollary \ref{cor:proper}) that the action of $G_{\mathbb{C}}$ on $M_{s}$ is
proper. (Keep in mind that we are assuming that $0$ is a regular value of the
moment map and that $G$ acts freely on $\Phi^{-1}(0).$) Because the action of
$G_{\mathbb{C}}$ is free, proper and holomorphic, the quotient
$M_{s}/G_{\mathbb{C}}$ has the structure of a complex manifold. On the other
hand, since each $G_{\mathbb{C}}$-orbit in $M_{s}$ intersects $\Phi^{-1}(0)$
in precisely one $G$-orbit, we have a natural bijective identification%
\begin{equation}
\Phi^{-1}(0)/G=M_{s}/G_{\mathbb{C}}.\label{identify}%
\end{equation}
Once we know that $M_{s}/G_{\mathbb{C}}$ has the structure of a complex
manifold, it is not hard to see that this complex structure agrees (under the
above identification) with the one obtained infinitesimally in
(\ref{infinitesimal complex}). (If $\pi_{\mathbb{C}}:M_{s}\rightarrow
M_{s}/G_{\mathbb{C}}$ is the quotient map and we identify $M_{s}%
/G_{\mathbb{C}}$ with $\Phi^{-1}(0)/G,$ then $\pi_{\mathbb{C}}$ agrees with
$\pi$ on $\Phi^{-1}(0). $ Thus $(\pi_{\mathbb{C}})_{\ast}$ agrees with
$\pi_{\ast}$ on vectors tangent to $\Phi^{-1}(0).$ But it is not hard to show
that every $(1,0)$-vector at a point in $\Phi^{-1}(0)$ is the sum of a
$(1,0)$-vector tangent to $\Phi^{-1}(0)$ and a $(1,0)$-vector tangent to the
$G_{\mathbb{C}} $ orbit.)

Although we do not require this result, we explain briefly how the quotient
$M_{s}/G_{\mathbb{C}}$ can be identified with the quotient in geometric
invariant theory (GIT). In GIT, there are several (in general, inequivalent)
notions of stable points for the action of a complex reductive group (such as
$G_{\mathbb{C}}$) on a compact K\"{a}hler manifold. (In \cite{Fogarty}, these
are called semistable, stable, and properly stable points.) Under our
assumptions---that $0$ is a regular value of the moment map and that $G$ acts
freely on $\Phi^{-1}(0)$---these different notions of stability in GIT turn
out to be equivalent to one another and to the notion of stability in
(\ref{eqn:stable def}). (As a consequence of a result of Kempf and Ness
\cite{Kempf-Ness}, the different notions of stability of $y\in M$ in GIT are
equivalent to: (1) the closure of $G_{\mathbb{C}}\cdot y$ intersects
$\Phi^{-1}(0),$ (2) $G_{\mathbb{C}}\cdot y$ itself intersects $\Phi^{-1}(0),$
and (3) $G_{\mathbb{C}}\cdot y$ intersects $\Phi^{-1}(0)$ and the stabilizer
of $y$ is finite. Condition 2 is our definition of the stable set. Since in
our case $M_{s}$ is open and acted on freely by $G_{\mathbb{C}},$ the other
two conditions are equivalent to Condition 2.) The GIT quotient coincides, as
a set and topologically, with $M_{s}/G_{\mathbb{C}}.$ See for example Theorem
8.3 of \cite{Fogarty} or Section 3 of \cite{Heinzner}.

The stable set can also be characterized in terms of the $G$-invariant
holomorphic sections of $\ell^{\otimes k}$ (see for example Donaldson's work
in \cite{Donaldson}, or \cite[Theorem 5.5]{Guillemin-Sternberg}):%
\begin{equation}
M_{s}=\{x\in M:s(x)\neq0\text{ for some }s\in\mathcal{H}(M;\ell^{\otimes
k})^{G}\text{ for some }k\}.\label{eqn:stable set character}%
\end{equation}
Suppose, then, that $s$ is an element of $s\in\mathcal{H}(M;\ell^{\otimes
k})^{G}$ that is not identically zero. (That such a section exists for some
$k$ follows from (\ref{eqn:stable set character}) and is established in
\cite[Appendix]{Guillemin-Sternberg}.) Then the unstable set $M_{s}\setminus
M$ is contained in the zero-set of $s$ and is therefore a set of codimension
at least one. This implies, for example, that the unstable set has measure zero.

We now establish a structure result for $M_{s}$ that implies, among other
things, that the action of $G_{\mathbb{C}}$ on $M_{s}$ is proper.\ Let
$\pi_{\mathbb{C}}:M_{s}\rightarrow M/\!\!/G$ be the complex quotient map under
the identification of $M_{s}/G_{\mathbb{C}}$ with $M/\!\!/G$. That is to say,
$\pi_{\mathbb{C}}(z)=\pi(x)$ if $x$ belongs to the unique $G$-orbit in
$\Phi^{-1}(0)$ contained in the $G_{\mathbb{C}}$-orbit of $z.$

\begin{theorem}
\label{thm:gc bundle}The map $\pi_{\mathbb{C}}:M_{s}\rightarrow M/\!\!/G$ is
smooth (in particular, continuous), and $M_{s}$ has the structure of a smooth
principal $G_{\mathbb{C}}$-bundle with respect to this map.
\end{theorem}

\begin{corollary}
\label{cor:proper}The action of $G_{\mathbb{C}}$ on $M_{s}$ is proper.
\end{corollary}

\begin{corollary}
\label{cor:lambda}The map $\Lambda:\mathfrak{g}\times\Phi^{-1}(0)\rightarrow
M_{s}$ given by%
\[
\Lambda(\xi,x)=e^{i\xi}\cdot x
\]
is bijective and a diffeomorphism.
\end{corollary}

\begin{proof}[Proof of Corollary \protect\ref{cor:proper}]
Suppose we have $y_{k}\in M_{s}$ and $g_{k}\in G_{\mathbb{C}}$ such that $%
y_{k}$ converges to $y\in M_{s}$ and such that $g_{k}\cdot y_{k}$ converges
to $z\in M_{s}.$ We will show that $g_{k}$ must then be convergent, which
implies that the action of $G_{\mathbb{C}}$ on $M_{s}$ is proper \cite[pp.
53 \textit{ff. }and Chap. 2]{Duistermaat-Kolk}. Note that $\pi_{\mathbb{C}%
}(g_{k}y_{k})=\pi_{\mathbb{C}}(y_{k}).$ Taking limits and using the
continuity of $\pi_{\mathbb{C}}$ gives that $\pi_{\mathbb{C}}(z)=\pi_{%
\mathbb{C}}(y).$ Choose a neighborhood $U$ of $\pi_{\mathbb{C}}(z)$ such
that $\pi_{\mathbb{C}}^{-1}(U)$ is homeomorphic to $G_{\mathbb{C}}\times U$
in such a way that the $G_{\mathbb{C}}$ action corresponds to left
multiplication in $G_{\mathbb{C}}.$ Under this homeomorphism, $y_{k}$
corresponds to $(h_{k},u_{k})$ and $y$ corresponds to $(h,u),$ with $%
h_{k}\rightarrow h$ and $u_{k}\rightarrow u.$ Then $g_{k}\cdot y_{k}$
corresponds to $(g_{k}h_{k},u_{k})$ and $z$ corresponds to $(f,u),$ with $%
g_{k}h_{k}\rightarrow f.$ Thus, $g_{k}$ converges to $h^{-1}f.$
\end{proof}

\begin{proof}[Proof of Corollary \protect\ref{cor:lambda}]
The bijectivity of $\Lambda$ follows from: (1) the freeness of the action of
$G_{\mathbb{C}}$ on $M_{s},$ (2) the fact that each $G_{\mathbb{C}}$-orbit
intersects $\Phi^{-1}(0)$ in a single $G$-orbit, and (3) the bijectivity of
the Cartan decomposition of $G_{\mathbb{C}}$ \cite[Thm 6.31, pp. 362]{Knapp}.
To see that $\Lambda$ is a diffeomorphism, decompose $\pi_{\mathbb{C}%
}^{-1}(U)$ as $G_{\mathbb{C}}\times U.$ Since the Cartan decomposition of $%
G_{\mathbb{C}}$ is a diffeomorphism, we obtain $\pi_{\mathbb{C}%
}^{-1}(U)\simeq\mathfrak{g}\times G\times U.$ The quotient map $\pi:\Phi
^{-1}(0)\rightarrow M/\!\!/G$ is a principal $G$-bundle, so there is a slice
$N$ in $\Phi^{-1}(0)$ which is transverse to the $G$-action and diffeomorphic
to $U$. The free $G$-action on $\Phi^{-1}(0)$ then yields a diffeomorphism
from $G\times N$ to a $G$-invariant neighborhood $M_{0}$ in the zero-set.
Hence we obtain a local diffeomorphism $\mathfrak{g}\times M_{0}\simeq%
\mathfrak{g}\times G\times N\simeq\pi_{\mathbb{C}}^{-1}(U)$ which
corresponds to the map $\left. \Lambda\right\vert _{\mathfrak{g}\times
M_{0}} $. This shows, in particular, that the Jacobian of $\Lambda$ is
invertible at each point and hence that $\Lambda$ is a global diffeomorphism.
\end{proof}

\begin{proof}[Proof of Theorem \protect\ref{thm:gc bundle}]
Given a point $u$ in $\Phi^{-1}(0)/G,$ choose a neighborhood $U$ of $u,$ a
point $x\in\Phi^{-1}(0)$ with $\pi(x)=u,$ and an embedded submanifold $N$ through
$x$ which is transverse to the $G$-action so that $\pi$ maps $N$ injectively
and diffeomorphically onto $U.$ Let $\pi^{-1}:U\rightarrow N$ be the smooth inverse
to $\pi.$ Now consider the map $\Psi:G_{\mathbb{C}}\times U\rightarrow M_{s}$ given by%
\begin{equation*}
\Psi(g,u)=g\cdot\pi^{-1}(u).
\end{equation*}

First, we establish that $\Psi$ is injective (globally, not just locally).
If $u_{1}$ and $u_{2}$ are distinct elements of $N,$ then $\pi^{-1}(u_{1})$
and $\pi^{-1}(u_{2})$ are in distinct $G$-orbits. But distinct $G$-orbits in
$\Phi^{-1}(0)$ lie in distinct $G_{\mathbb{C}}$-orbits (see the comments
following the proof of Theorem \ref{thm:s along a curve}(a)). Thus, $%
\Psi(g_{1},u_{1})\neq\Psi(g_{2},u_{2})$ if $u_{1}\neq u_{2}.$ Then if $%
g_{1}\neq g_{2},$ $\Psi(g_{1},u)\neq\Psi (g_{2},u),$ because $G_{\mathbb{C}}$
acts freely on $M_{s}.$

Next, we establish that the differential of $\Psi$ is injective at each
point. Since $\pi$ maps $N$ diffeomorphically onto $U,$ it suffices to prove
the same thing for the map $\mathrm{\Gamma}:G_{\mathbb{C}}\times
N\rightarrow M_{s}$ given by $\mathrm{\Gamma}(g,n)=g\cdot n.$ So consider $%
(g,n)\in G_{\mathbb{C}}\times N$. Since $\pi$ maps $N$ diffeomorphically
onto $U,$ $T_{n}(\Phi^{-1}(0))$ is the direct sum of $T_{n}(N)$ and $%
T_{n}(G\cdot n).$ Furthermore, for any $\xi \in\mathfrak{g},$ the vector
field $JX^{\xi}$ is nonzero at $n$ and orthogonal to the tangent space to $%
\Phi^{-1}(0).$ It follows that%
\begin{equation*}
T_{n}(M)=T_{n}(G_{\mathbb{C}}\cdot n)\oplus T_{n}(N)
\end{equation*}
and thus%
\begin{align*}
T_{gn}(M) & =g_{\ast}(T_{n}(G_{\mathbb{C}}\cdot n))\oplus g_{\ast}(T_{n}(N))
\\
& =T_{gn}(G_{\mathbb{C}}\cdot n)\oplus g_{\ast}(T_{n}(N))
\end{align*}
since the action of $g$ takes $G_{\mathbb{C}}\cdot n$ into itself.
Now suppose that $\alpha(t)=(g(t),n(t))$ is a curve in $G_{\mathbb{C}}\times
N$ passing through $(g,n)$ at $t=0$ and such that $g^{\prime}(0)g^{-1}=\xi
\in\mathfrak{g}_{\mathbb{C}}$ and $n^{\prime}(t)=X.$ Then applying $\mathrm{%
\Gamma}$ and differentiating at $t=0$ gives
\begin{equation*}
\mathrm{\Gamma}_{\ast}(\alpha^{\prime}(0))=X^{\xi}_{g\cdot n}+g_{\ast}(X),
\end{equation*}
which is nonzero provided that $\alpha^{\prime}(0)$ is nonzero.

Finally, we prove the theorem. Since $\Psi$ is globally injective, it has an
inverse map defined on its image, where the image of $\Psi$ is simply $\pi_{%
\mathbb{C}}^{-1}(U)$. The inverse function theorem tells us, then, that $%
\pi_{\mathbb{C}}^{-1}(U)$ is an open set and that the inverse map to $\Psi$
is smooth. Thus, $\pi_{\mathbb{C}}^{-1}(U)$ is diffeomorphic to $G_{\mathbb{C%
}}\times U$ in such a way that the action of $G_{\mathbb{C}}$ corresponds to
the left action of $G_{\mathbb{C}}$ on itself. Under this diffeomorphism, $%
\pi_{\mathbb{C}}$ corresponds to projection onto the second factor, which is
smooth. Thus the diffeomorphism $\Psi$ shows that $M_{s}\rightarrow M/\!\!/G$
has the necessary smooth local triviality property to be a smooth principal $%
G_{\mathbb{C}}$-bundle.
\end{proof}

\subsection{Quantum reduction and the Guillemin--Sternberg map}

In the previous section, we considered reduction at the classical level, which
amounts to passing from $M$ to $M/\!\!/G.$ Alternatively, we may first
quantize $M$, by looking at the space of holomorphic sections of $l^{\otimes
k}$ over $M,$ and then perform reduction at the quantum level. According to
the philosophy of Dirac \cite[Lecture 2, pp. 34 \emph{ff}.]{Dirac}, reduction at the quantum level amounts
to looking at holomorphic sections that are invariant under an appropriate
action of the group $G,$ that is, the quantum reduced space is the null-space
of the quantized moment map ((\ref{eqn:G action on sections}) below).

We now describe how the action of $G$ on the space of sections is constructed.
Following the program of geometric quantization, we first define an action of
the Lie algebra $\mathfrak{g}$ on the space of smooth sections of $\ell$ by
\[
Q_{\xi}:=\nabla_{X^{\xi}}-i\,\phi_{\xi},\quad\xi\in\mathfrak{g}.
\]
These operators satisfy $[Q_{\xi},Q_{\eta}]=Q_{[\xi,\eta]}$. The prequantum
bundle $\ell$ is said to be $G$-invariant if this action of $\mathfrak{g}$ can
be exponentiated to an action of the group $G.$ (If $G$ is simply connected,
this can always be done.) We henceforth assume that $\ell$ is $G$-invariant.
Since $\ell$ is a holomorphic line bundle (i.e.,the total space is a complex
manifold), it is not hard to show that the $G$-action on $\ell$ can be
analytically continued to an action of $G_{\mathbb{C}}$ (by following an
argument similar to that of Guillemin and Sternberg that the action of $G$ on
$M$ can be analytically continued to an action of $G_{\mathbb{C}}$ on $M$
\cite[Theorem 4.4]{Guillemin-Sternberg}).

For each $k,$ we then define an action of $\mathfrak{g}$ on the space of
smooth sections of $\ell^{\otimes k}$ by%
\begin{equation}
Q_{\xi}:=\nabla_{X^{\xi}}-ik\phi_{\xi},\quad\xi\in\mathfrak{g}%
.\label{eqn:G action on sections}%
\end{equation}
(We suppress the dependence on $k.$) These also satisfy $[Q_{\xi},Q_{\eta
}]=Q_{[\xi,\eta]}$ and they preserve the quantum Hilbert space $\mathcal{H}%
(M;\ell^{\otimes k})$ of \emph{holomorphic} sections of $\ell^{\otimes k}.$
Since $\ell$ is $G$-invariant, it follows that these operators can be
exponentiated to an action of $G$ that preserves $\mathcal{H}(M;\ell^{\otimes
k}).$ The space obtained by performing reduction at the quantum level is then
the space of $G$-invariant holomorphic sections of $\ell^{\otimes k},$ denoted
$\mathcal{H}(M;\ell^{\otimes k})^{G}.$

We see, then, that if we first quantize $M$ and then reduce by $G,$ we obtain
the space $\mathcal{H}(M;\ell^{\otimes k})^{G}$ of $G$-invariant sections of
$\ell^{\otimes k}$ over $M.$ On the other hand, we may first reduce $M$ by $G$
at the classical level and then perform quantization of the reduced manifold
$M/\!\!/G.$ Assuming that the bundle $\ell$ is $G$-invariant, it is not hard
to see that it descends naturally to a Hermitian line bundle $\hat{\ell}$ with
connection over $M/\!\!/G,$ whose curvature is equal to $-i\widehat{\omega}.$
The bundle $\hat{\ell}$ can be made into a holomorphic line bundle in the same
way as $\hat{\ell},$ by decreeing that the holomorphic sections are those that
are covariantly constant in the $(0,1)$-directions. The space, then, obtained
by first reducing and then quantizing is the space of all holomorphic sections
of $\hat{\ell}^{\otimes k}$ over $M/\!\!/G,$ denoted $\mathcal{H}%
(M/\!\!/G;\hat{\ell}^{\otimes k}).$

In the paper \cite{Guillemin-Sternberg}, Guillemin and Sternberg consider a
geometrically natural linear map $A_{k}$ between the \textquotedblleft first
quantize and then reduce\textquotedblright\ space $\mathcal{H}(M;\ell^{\otimes
k})^{G}$ and the \textquotedblleft first reduce and then
quantize\textquotedblright\ space $\mathcal{H}(M/\!\!/G;\hat{\ell}^{\otimes
k}).$ This map consists simply of taking a $G$-invariant holomorphic section
over $M,$ restricting it to $\Phi^{-1}(0),$ and then letting it descend to
$M/\!\!/G=\Phi^{-1}(0)/G.$ The remarkable result established in \cite[Thm
5.2]{Guillemin-Sternberg} is that this natural map is invertible. This result
has been generalized in various ways, to symplectic manifolds that are not
K\"{a}hler and to situations where the quotient is singular. See, for example,
\cite{Huebschmann}, \cite{Jeffrey}, \cite{Sjamaar1}, \cite{Sjamaar2},
\cite{Meinrenken}, and \cite{Tian-Zhang}.

In view of the importance of the map $A_{k}$, we briefly sketch the proof of
its invertibility. It is easy to see that a section of $\hat{\ell}$\ defines a
$G$-invariant section of $\ell$ over the zero-set by pullback along the
canonical projection $\pi:\Phi^{-1}(0)\rightarrow M/\!\!/G$, and vice-versa.
Moreover, a $G$-invariant section over the zero-set has a unique analytic
continuation to the stable set, obtained from the $G_{\mathbb{C}}$-action on
$M$. The difficulty is in showing that the resulting analytic continuation
extends smoothly across the unstable set.\ To construct this extension,
Guillemin and Sternberg show that the magnitude of a $G$-invariant holomorphic
section over the stable set approaches zero as one approaches the unstable
set. The Riemann Extension Theorem then shows that such a section extends
holomorphically across the unstable set.

The invertibility of the map $A_{k}$ shows that, \textit{in a certain sense},
quantization commutes with reduction. That is, $A_{k}$ gives a natural
identification of the \textit{vector space} obtained by first quantizing and
then reducing with the \textit{vector space} obtained by first reducing and
then quantizing. However, in quantum mechanics, the structure of the spaces as
\textit{Hilbert} spaces, not just vector spaces, is essential. For example,
the expectation value of an operator $D$ in the \textquotedblleft
state\textquotedblright\ $s$ is given by $\left\langle s,Ds\right\rangle ,$
where $s$ is a unit vector in the relevant Hilbert space. It is natural, then,
to consider the extent to which the map $A_{k}$ is unitary with respect to the
natural inner products (\ref{innerproduct}) and (\ref{innerproduct2}) on the
two spaces.

Our Theorem \ref{thm:s norm decomp} suggests that it is very rare for $A_{k}$
(or any constant multiple of $A_{k}$) to be unitary. If the $G$-orbits in
$\Phi^{-1}(0)$ have nonconstant volume, Theorem \ref{thm:non-unitarity} shows
that for all sufficiently large $k,$ $A_{k}$ is not a constant multiple of a
unitary map. Indeed, in such cases, $A_{k}$ is not even asymptotic to a
multiple of a unitary map in the limit as Planck's constant tends to zero
(i.e., as $k$ tends to infinity). We may say therefore that, in the setting of
\cite{Guillemin-Sternberg}, quantization does \textit{not} commute with
reduction in the strongest desirable sense.

\section{\label{subsec:metaplectic correction}A map of Guillemin--Sternberg
type in the presence of the metaplectic correction}

In this section, we consider the \textquotedblleft metaplectic
correction,\textquotedblright\ which consists of tensoring the original line
bundle $\ell$ by the square root of the canonical bundle of $M,$ and similarly
for $\hat{\ell}.$ We introduce here an analog $B_{k}$ of the
Guillemin--Sternberg map in the presence of the metaplectic correction and we
show (Theorem \ref{thm:modGS invertible} of this section) that $B_{k}$ is
invertible for all sufficiently large $k.$ We will see eventually (Theorems
\ref{thm:Bk unitary} and \ref{thm:non-unitarity}) that the maps $B_{k},$
unlike the maps $A_{k},$ become approximately unitary as Planck's constant
tends to zero (i.e., as $k$ tends to infinity).

There are two conceptual differences between the original map $A_{k}$ and the
\textquotedblleft metaplectically corrected\textquotedblright\ map $B_{k}.$
First, $G$-invariant holomorphic sections of the square root of the canonical
bundle behave badly as one approaches the unstable set. To compensate for
this, we need rapid decay of holomorphic sections of $G$-invariant holomorphic
sections of $\ell^{\otimes k}.$ This means that we obtain invertibility of the
maps $B_{k}$ for all sufficiently large $k,$ rather than for all $k$, as is
the case for the maps $A_{k}.$ Second, the pointwise magnitude of a section
$B_{k}s$ does not agree with the pointwise magnitude of the original section
$s$ on $\Phi^{-1}(0),$ in contrast to the map $A_{k}.$ Rather, the pointwise
magnitudes differ by a factor involving the volume of the $G$-orbits in
$\Phi^{-1}(0)$ (Theorem \ref{thm:BKS comparison}). The volume factor in
Theorem \ref{thm:BKS comparison} ultimately
cancels an unwanted volume factor in the asymptotics of the maps $A_{k}$
(Section \ref{sec:asymptotics}), allowing the maps $B_{k}$ to be asymptotically
unitary even though maps $A_{k}$ are not.

\subsection{Half-form bundles on $M$ and $M/\!\!/G$}

Let \thinspace$K=\bigwedge\nolimits^{n}\left(  T^{1,0}M\right)  ^{\ast}$
denote the canonical bundle of $M,$ that is, the top exterior power of the
bundle of $(1,0)$-forms. A smooth section of $K$ is called an $(n,0)$-form,
and the set of such forms is denoted $\Omega^{n,0}(M).$ An $(n,0)$-form is
called holomorphic if in each holomorphic local coordinate system, the
coefficient of $dz_{1}\wedge\cdots\wedge dz_{n}$ is a holomorphic function.
Equivalently, we may define a \textquotedblleft partial
connection\textquotedblright\ (defined only for vector fields of type $(0,1)$)
on the space of $(n,0)$-forms by setting%
\[
\nabla_{X}\alpha=i_{X}(d\alpha)
\]
whenever $X$ is of type $(0,1).$ It is easily verified that $\alpha$ is
holomorphic if and only if $\nabla_{X}\alpha=0$ for all vector fields of type
$(0,1).$ (Compare \cite[Sect. 9.3]{Woodhouse}.)

A choice of square root $\sqrt{K}$ of the canonical bundle, if it exists, is
called a half-form bundle.\ Since the first Chern class of the canonical
bundle is $-c_{1}(M)$, the canonical bundle will admit a square root if and
only if $-c_{1}(M)/2$ is an integral class; that is, if and only if the second
Stiefel--Whitney class $w_{2}(M)$ (which is the reduction $\operatorname{mod}$
$2$ of the first Chern class for a complex manifold) vanishes.\ We now assume
that $-c_{1}(M)/2$ is integral and we fix a choice of $\sqrt{K}$. (It is
likely that results similar to the ones in this paper could be obtained
assuming that $[\omega/2\pi]-c_{1}(M)/2$ is integral, rather than assuming, as
we do, that $[\omega/2\pi]$ and $c_{1}(M)/2$ are separately integral. See
\cite{Czyz} or \cite[Sect. 10.4]{Woodhouse}.)

One can define a partial connection acting on sections $\nu$ of $\sqrt{K}$ by
requiring that%
\begin{equation}
2\left(  \nabla_{X}\nu\right)  \nu=\nabla_{X}(\nu^{2}).\label{partial half}%
\end{equation}
(See again \cite[Sect. 9.3]{Woodhouse}.) The bundle $\sqrt{K}$ can then be
made into a holomorphic line bundle by defining the holomorphic sections of
$\sqrt{K}$ to be those for which $\nabla_{X}\nu=0$ for all vector fields of
type $(0,1).$

Because the action of $G$ on $M$ is holomorphic, the action of $G$ on
$n$-forms preserves the space of $(n,0)$-forms and the space of holomorphic
$n$-forms. The associated action of a Lie algebra element $\xi\in\mathfrak{g}$
on $(n,0)$-forms is by the Lie derivative $\mathcal{L}_{X^{\xi}}.$ There is an
associated action of $\mathfrak{g}$ on the space of sections of $\sqrt{K}$,
also denoted $\mathcal{L}_{X^{\xi}}$, satisfying (by analogy with
(\ref{partial half}))%
\begin{equation}
2\left(  \mathcal{L}_{X^{\xi}}\nu\right)  \nu=\mathcal{L}_{X^{\xi}}(\nu
^{2}),\label{lie action}%
\end{equation}
and this action preserves the space of holomorphic sections.

There is a natural way to define a Hermitian structure on $\sqrt{K},$ which is
standard in geometric quantization. (It is a special case of the BKS pairing;
see \cite[Sec. 10.4]{Woodhouse}.) If $\nu,\mu\in\mathrm{\Gamma}(\sqrt{K}) $
are half-forms, then $\nu^{2}\wedge\overline{\mu}^{2}\in\mathrm{\Gamma
}(\bigwedge\nolimits^{2n}T_{\mathbb{C}}^{\ast}(M)).$ The volume form
$\varepsilon_{\omega}=\omega^{\wedge n}/n!$ is a global trivializing section
of the determinant bundle $\bigwedge\nolimits^{2n}T^{\ast}M$, and so there is
a function, denoted by $(\nu,\mu),$ such that
\begin{equation}
\nu^{2}\wedge\overline{\mu}^{2}=\left(  \nu,\mu\right)  ^{2}\,\varepsilon
_{\omega}.\label{eqn:Herm form on half-forms}%
\end{equation}

We can use this pairing to define a Hermitian form on the tensor product
$\ell^{\otimes k}\otimes\sqrt{K}$: \ for sections $t_{1},t_{2}\in
\mathrm{\Gamma} (\ell^{\otimes k}\otimes\sqrt{K})$ which are locally
represented by $t_{j}(x)=s_{j}(x)\mu(x)$, we define%
\[
(t_{1},t_{2})(x)=(s_{1}\nu,s_{2}\mu)(x)=(s_{1}(x),s_{2}(x))\left(  \nu
,\mu\right)  (x).
\]
Denote the pairing of a half-form with itself by $\left\vert \nu\right\vert
^{2}=\left(  \nu,\nu\right)  $.

Let $\widehat{K}$ denote the canonical bundle over the reduced manifold
$M/\!\!/G.$ We now assume that the action (\ref{lie action}) of $\mathfrak{g}$
on sections of $\sqrt{K}$ exponentiates to an action of the group $G$ (this is
automatic, for example, if $G$ is simply connected). In the next subsection,
we will show that this assumption allows us to construct in a natural way a
square root $\sqrt{\widehat{K}}$ of $\widehat{K}$ that is related in a nice
way to the chosen square root $\sqrt{K}$ of $K.$ Once $\sqrt{\widehat{K}}$ is
constructed, we will define a Hermitian structure on $\sqrt{\widehat{K}}$
analogously to (\ref{eqn:Herm form on half-forms}):
\begin{equation}
\nu^{2}\wedge\overline{\mu}^{2}=\left(  \nu,\mu\right)  ^{2}\,\varepsilon
_{\widehat{\omega}},\label{eqn:half-form Hermitian form}%
\end{equation}
where $\varepsilon_{\widehat{\omega}}=\widehat{\omega}^{\wedge(n-d)}/(n-d)!$
is the volume form on $M/\!\!/G.$

\subsection{The modified Guillemin--Sternberg
map\label{subsec:modified GS map}}

We continue to assume that the canonical bundle $K$ of $M$ admits a square
root and that we have chosen one such square root and denoted it by $\sqrt
{K}.$ We also continue to assume that the action of the Lie algebra
$\mathfrak{g}$ on sections of $\sqrt{K},$ given by (\ref{lie action}),
exponentiates to an action of the group $G.$

In this subsection, we describe a natural map $B_{k}$ from the space of
$G$-invariant holomorphic sections of $\ell^{\otimes k}\otimes\sqrt{K}$ to the
space of holomorphic sections of $\hat{\ell}^{\otimes k}\otimes\sqrt
{\widehat{K}}.$ This map is the analog of the Guillemin--Sternberg map $A_{k}$
between the $G$-invariant holomorphic sections of $\ell^{\otimes k}$ and the
space of holomorphic sections of $\hat{\ell}^{\otimes k}.$ However, the
construction of $B_{k}$ is more involved than that of $A_{k}.$ After all,
sections of $\sqrt{K}$ are square roots of $(n,0)$-forms, whereas sections of
$\sqrt{\widehat{K}}$ are square roots of $(n-d,0)$-forms. Thus the map $B_{k}$
must include a mechanism for changing the degree of a half-form.

There is one further, vitally important difference between the map $B_{k}$ and
the map $A_{k}.$ The bundle $\hat{\ell}$ inherits its Hermitian structure from
that of $\ell.$ As a result, the pointwise magnitude of a section a
$G$-invariant section of $\ell^{\otimes k}$ is the same on $\Phi^{-1}(0)$ as
the pointwise magnitude of the corresponding section of $\hat{\ell}^{\otimes
k}.$ That is, for each $s\in\mathcal{H}(M;\ell^{\otimes k})^{G}$ and each
$x_{0}\in\Phi^{-1}(0)$ we have%
\begin{equation}
\left\vert s\right\vert ^{2}(x_{0})=\left\vert A_{k}s\right\vert ^{2}%
([x_{0}]).\label{magnitude1}%
\end{equation}

By contrast, $\sqrt{\widehat{K}}$ has its own intrinsically defined Hermitian
structure given by (\ref{eqn:Herm form on half-forms}) It turns out, then,
that $B_{k}$ does not satisfy the analog of (\ref{magnitude1}). Rather, we
will show (Theorem \ref{thm:BKS comparison}) that for $r\in\mathcal{H}%
(M;\ell^{\otimes k}\otimes\sqrt{K})^{G}$ and each $x_{0}\in\Phi^{-1}(0)$ we
have%
\begin{equation}
\left\vert r\right\vert ^{2}(x_{0})=2^{d/2}\mathrm{\operatorname*{vol}}(G\cdot
x_{0})^{-1}\left\vert B_{k}r\right\vert ^{2}([x_{0}]).\label{magnitude2}%
\end{equation}
We will see eventually that the volume factor in (\ref{magnitude2}) cancels an
unwanted volume factor in the asymptotic behavior of $A_{k}.$ The result is
that maps $B_{k}$ (unlike the maps $A_{k}$) are asymptotically unitary as $k$
tends to infinity.

Finally, we address the invertibility of the map $B_{k}.$ Again, the situation
is slightly different from that with $A_{k}.$ Half-forms typically have bad
behavior near the unstable set, and as a result, we are only able to prove
that $B_{k}$ is invertible for \textit{sufficiently large }$k.$ (This is in
contrast to the map $A_{k}$, which is invertible for all $k$.)

\bigskip

Before we get to half-forms, it will be useful to describe the descent of
$G_{\mathbb{C}}$-invariant $(n,0)$-forms on the stable set $M_{s}$ to $(n-d,0)
$-forms on $M/\!\!/G.$ We cannot simply restrict an $(n,0)$-form to the
zero-set and then let it descend to $M/\!\!/G,$ as we did for the bundle
$\ell^{\otimes k},$ because the result would not be an $(n-d)$-form on
$M/\!\!/G$. The process, which we describe in detail below, is to first
contract with the infinitesimal $G$-directions and then use $\pi_{\mathbb{C}}$
to push the result down to the quotient. It turns out that this process is
invertible and preserves holomorphicity.

Choose an $\operatorname*{Ad}$-invariant inner product on $\mathfrak{g}$ which
is normalized so that the volume of $G$ with respect to the associated Haar
measure is $1$. Then fix a basis $\mathrm{\Xi}=\{\xi_{1},...,\xi_{d}\}$ of the
Lie algebra $\mathfrak{g}$ which is orthonormal with respect to this inner
product. Given a $G_{\mathbb{C}}$-invariant $(n,0)$-form defined on the stable
set $M_{s}$, define an $(n-d,0)$-form $\beta$ on $M_{s}$ by%
\[
\beta=\mathfrak{i}(\bigwedge\nolimits_{j}X^{\xi_{j}})\alpha.
\]
We claim that $\beta$ has the properties that for each $X\in T(G_{\mathbb{C}%
}\cdot x)$%
\begin{align}
\,\mathfrak{i}(X)\beta &  =0~\text{and }%
\label{eqn:(n,0)-form reduction props1}\\
\mathfrak{i}(X)d\beta &  =0.\label{eqn:(n,0)-form reduction props2}%
\end{align}
It is clear that (\ref{eqn:(n,0)-form reduction props1}) holds when $X$ is of
the form $X^{\xi}$ with $\xi\in\mathfrak{g}.$ But since $\beta$ is an
$(n-d,0)$-form, $\mathfrak{i}(JX^{\xi})\beta=\sqrt{-1}\mathfrak{i}(X^{\xi
})\beta,$ and so (\ref{eqn:(n,0)-form reduction props1}) holds for $X=X^{\xi}
$ with $\xi\in\mathfrak{g}_{\mathbb{C}}.$ To verify
(\ref{eqn:(n,0)-form reduction props2}), we first note that in the presence of
(\ref{eqn:(n,0)-form reduction props1}),
(\ref{eqn:(n,0)-form reduction props2}) is equivalent to $\mathcal{L}%
_{X}(\beta)=0,$ i.e., to the condition that $\beta$ be $G_{\mathbb{C}}%
$-invariant. So we need to verify that $\beta$ is invariant if $\alpha$ is.
Since $\alpha$ is of type $(n,0),$ contracting $\alpha$ with
$\bigwedge\nolimits_{j}X^{\xi_{j}}$ is the same as contracting it with
$\bigwedge\nolimits_{j}\pi_{+}X^{\xi_{j}},$ and a simple computation shows
that the polyvector $\bigwedge\nolimits_{j}\pi_{+}X^{\xi_{j}}$ is
$G_{\mathbb{C}}$-invariant. Contracting a $G_{\mathbb{C}}$-invariant form with
a $G_{\mathbb{C}}$-invariant polyvector gives another $G_{\mathbb{C}}%
$-invariant form.

Now, given any $(n-d,0)$-form $\beta$ on $M_{s}$ satisfying
(\ref{eqn:(n,0)-form reduction props1}) and
(\ref{eqn:(n,0)-form reduction props2}), it is not hard to show that there is
a unique $(n-d,0)$-form $\widehat{\beta}$ on $M/\!\!/G=M_{s}/G_{\mathbb{C}}$
such that $\pi_{\mathbb{C}}^{\ast}(\widehat{\beta})=\beta$, where
$\pi_{\mathbb{C}}^{\ast}:M_{s}\rightarrow M_{s}/G_{\mathbb{C}}$ is the
quotient map. So we have given a procedure for turning a $G_{\mathbb{C}}%
$-invariant $(n,0)$-form $\alpha$ on $M$ into an $(n-d,0)$-form $\widehat
{\beta}$ on $M/\!\!/G.$

In the other direction, suppose $\widehat{\beta}$ is an $(n-d,0)$-form on
$M/\!\!/G.$ Then the pullback $\beta:=\pi_{\mathbb{C}}^{\ast}(\widehat{\beta
})$ is an $(n-d,0)$-form on $M_{s}$ that (as is easily verified) satisfies
(\ref{eqn:(n,0)-form reduction props1}) and
(\ref{eqn:(n,0)-form reduction props2}). We can construct from $\beta$ an
$(n,0)$-form $\alpha$ as follows: given a local frame $\{X^{\xi_{1}}%
,\dots,X^{\xi_{d}},Y_{1},\dots,Y_{n-d}\}$ for $T_{x}M_{s}$, set
\begin{equation}
\alpha(X^{\xi_{1}},\dots,X^{\xi_{d}},Y_{1},\dots,Y_{n-d})=\pi_{\mathbb{C}%
}^{\ast}\widehat{\beta}(Y_{1},\dots,Y_{n-d})\label{eqn:(n-d,0)-form lift}%
\end{equation}
and define $\alpha$ on any other frame by $GL(n,\mathbb{C)}$-equivariance and
the requirement that $\alpha$ be an $(n,0)$-form. (Note that the tangent space
at a point in the stable set is a direct sum of the tangent space to the
$G_{\mathbb{C}}$-orbit through that point and the transverse directions, so
that every frame is equivalent to a linear combination of frames which are
$GL(n,\mathbb{C)}$-equivalent to one of the form $(W_{1},\dots,W_{d}%
,Y_{1},\dots,Y_{n-d})$ where $W_{j}=X^{\xi_{j}}$ or $JX^{\xi_{j}}$; we define
$\,\mathfrak{i}(JX^{\xi})\alpha=\sqrt{-1}\,\mathfrak{i}(X^{\xi})\alpha$.) It
is easily verified that $\alpha$ is again $G_{\mathbb{C}}$-invariant.

The two processes we have defined---contracting a $G_{\mathbb{C}}$-invariant
$(n,0)$-form on $M_{s}$ with $\bigwedge\nolimits_{j}X^{\xi_{j}}$ and the
letting it descent to the quotient, and pulling back an $(n-d,0)$-form on
$M/\!\!/G$ and then \textquotedblleft expanding\textquotedblright\ it by
(\ref{eqn:(n-d,0)-form lift})---are clearly inverse to each other. They
therefore define a bijective map
\begin{equation}
\alpha\in\Omega^{n,0}(M_{s})^{G_{\mathbb{C}}}\mapsto\mathfrak{B}(\alpha
)\in\Omega^{n-d,0}(M/\!\!/G)\label{eqn:restrict-then-contract map}%
\end{equation}
where $\mathfrak{B}(\alpha)$ is the unique $(n-d,0)$-form on $M/\!\!/G$ such
that $\pi_{\mathbb{C}}^{\ast}\mathfrak{B}(\alpha)=\mathfrak{i}(\bigwedge
\nolimits_{j}X^{\xi_{j}})\alpha.$ This bijective map has the further property
that $\alpha$ is locally holomorphic if and only if $\mathfrak{B}(\alpha)$ is
locally holomorphic since contracting with $\bigwedge\nolimits_{j}X^{\xi_{j}}$
is the same as contracting with $\bigwedge\nolimits_{j}\pi_{+}X^{\xi_{j}}$ and
the vector fields $\pi_{+}X^{\xi}$ are holomorphic.

Finally, note that the contraction $\mathfrak{i}\left(  X^{\xi_{j}}\right)
\alpha$ is proportional to the contraction $\mathfrak{i}\left(  JX^{\xi_{j}%
}\right)  \alpha$. By the results of Section \ref{subsec:Kahler reduction},
pushing down a $G_{\mathbb{C}}$-invariant $(n,0)$-form by the above process is
equivalent to first restricting it to $\Phi^{-1}(0)$ and then pushing it down
by the analogous process using the quotient map $\pi:\Phi^{-1}(0)\rightarrow
M/\!\!/G\simeq\Phi^{-1}(0)/G$. The vectors $JX^{\xi}$ span the normal bundle
of $\Phi^{-1}(0)$ in $M$, so the contraction step can then be understood as
contracting with the directions normal to the zero-set; this is perhaps the
obvious way to \textquotedblleft restrict\textquotedblright\ a top-dimensional
form to a submanifold. Since the Guillemin--Sternberg map $A_{k}$ is defined
as \textquotedblleft restrict to $\Phi^{-1}(0)$ and then descend to
$M/\!\!/G$\textquotedblright, we will sometimes interpret the map
$\mathfrak{B}$ as first contracting, then restricting the result to $\Phi
^{-1}(0)$, and finally pushing the result to the quotient.

\bigskip

We now turn to the descent map for half-forms. Recall that we assume that the
infinitesimal action of $\mathfrak{g}$ on $\sqrt{K}$ defined by
(\ref{lie action}) exponentiates to an action of $G$ on $\sqrt{K}$ which is
compatible with the action of $G$ on $K$. Following essentially the same
argument that Guillemin and Sternberg make for the analytic continuation of
the action of $G$ on $M$, it can be shown that this action lifts to an action
of $G_{\mathbb{C}}$ which covers the $G_{\mathbb{C}}$-action on $K$.

A moment's thought shows that the map $\mathfrak{B}$ actually provides an
identification of $K_{x}$ with $\widehat{K}_{[x]},$ $x\in M_{s},$ and this
identification commutes with the action of $G_{\mathbb{C}},$ that is, for each
$g\in G_{\mathbb{C}}$%
\[
g^{\ast}\mathfrak{B}(\alpha)=\mathfrak{B}(g^{\ast}\alpha).
\]
The contraction map $\mathfrak{B}$ therefore identifies $\left.  K\right\vert
_{G_{\mathbb{C}}\cdot x}$ with $\widehat{K}_{[x]}.$

Let $\sqrt{\widehat{K}}$ denote a line bundle over $M/\!\!/G$ whose fiber is
the equivalence class of $\left.  \sqrt{K}\right\vert _{G_{\mathbb{C}}\cdot x}
$ under the $G_{\mathbb{C}}$-action. Since tensoring commutes with the
$G_{\mathbb{C}}$-action, it follows that this bundle is a square root of the
canonical bundle on $M/\!\!/G$. We let $\mathrm{\Gamma}(M,\sqrt{K})$ and
$\mathrm{\Gamma}(M/\!\!/G,\sqrt{\widehat{K}})$ denote the space of smooth
sections of $\sqrt{K}$ and $\sqrt{\widehat{K}},$ respectively, and we let
$\mathrm{\Gamma}(M,\sqrt{K})^{G}$ denote the space of $G$-invariant sections
of $\sqrt{K}.$

By construction, then, the pullback of the half-form bundle on the quotient by
the quotient map is isomorphic to the half-form bundle on the stable set;
i.e.,
\[
\pi_{\mathbb{C}}^{\ast}\sqrt{\widehat{K}}\simeq\left.  \sqrt{K}\right\vert
_{M_{s}}.
\]
Moreover, this isomorphism defines a map
\[
B:\mathrm{\Gamma}(M,\sqrt{K} )^{G_{\mathbb{C}}}\rightarrow\mathrm{\Gamma
}(M/\!\!/G,\sqrt{\widehat{K}})
\]
such that $(B\nu_{x},B\nu_{x})=\mathfrak{B}(\nu^{2}).$

We have thus proved the following.

\begin{theorem}
\label{thm:modified GS map}There exists a linear map $B:\mathrm{\Gamma
}(M,\sqrt{K})^{G}\rightarrow\mathrm{\Gamma}(M/\!\!/G,\sqrt{\widehat{K}}),$
unique up to an overall sign, with the property that%
\[
\pi_{\mathbb{C}}^{\ast}\left[  (B\nu)^{2}\right]  =\left.  \left[
\mathfrak{i}\left(  \bigwedge\nolimits_{j}X^{\xi_{j}}\right)  (\nu
^{2})\right]  \right\vert _{M_{s}}.
\]
For any open set $U$ in $M/\!\!/G,$ if $\nu$ is holomorphic in a neighborhood
$V$ of $\pi_{\mathbb{C}}^{-1}(U),$ then $B\nu$ is holomorphic on $U.$

For each $k,$ there is a linear map $B_{k}:\mathrm{\Gamma}(M,\ell^{\otimes
k}\otimes\sqrt{K})^{G}\rightarrow\mathrm{\Gamma}(M/\!\!/G,\hat{\ell}^{\otimes k}%
\otimes\sqrt{\widehat{K}}),$ unique up to an overall sign, with the property
that%
\[
B_{k}(s\otimes\nu)=A_{k}(s)\otimes B(\nu)
\]
for all $s\in\mathrm{\Gamma}(\ell^{\otimes k})$ and $\nu\in\mathrm{\Gamma
}(\sqrt{K}).$ This map takes holomorphic sections of $\left.  \ell^{\otimes
k}\otimes\sqrt{K}\right\vert _{V}$ to holomorphic sections of $\left.
\hat{\ell}^{\otimes k}\otimes\sqrt{\widehat{K}}\right\vert _{U}.$
\end{theorem}

The last claim in Theorem \ref{thm:modified GS map} follows from the
definition (\ref{partial half}) of the partial connection on $\sqrt{K}.$

We obtain a version of Guillemin and Sternberg's \textquotedblleft
quantization commutes with reduction\textquotedblright using the modified map
$B_{k}$. Our version is necessarily weaker since the pointwise behavior of
$G$-invariant elements of $\mathcal{H}(M;\ell^{\otimes k}\otimes\sqrt{K})$ is
worse that the behavior $G$-invariant sections of $\mathcal{H}(M;\ell^{\otimes
k})$ (Theorem \ref{thm:s along a curve}).

\begin{theorem}
\label{thm:modGS invertible}For $k$ sufficiently large, the map
\[
B_{k}:\mathcal{H}(M;\ell^{\otimes k}\otimes\sqrt{K})^{G}\rightarrow
\mathcal{H}(M/\!\!/G;\hat{\ell}^{\otimes k}\otimes\sqrt{\widehat{K}})
\]
is bijective.
\end{theorem}

The proof of this result (below) is similar to the proof in
\cite{Guillemin-Sternberg} of the invertibility of the map $A_{k},$ with a few
modifications to deal with the half-forms.

The inclusion of half-forms modifies the Hermitian form on sections of the
reduced corrected prequantum bundle $\hat{\ell}^{\otimes k}\otimes
\sqrt{\widehat{K}}$ in a way that is not proportional to the modification of
Hermitian form on sections of $\ell^{\otimes k}\otimes\sqrt{K}$. We can
actually compute the difference:

\begin{theorem}
\label{thm:BKS comparison}Suppose $r\in\mathcal{H}(M;\ell^{\otimes k}%
\otimes\sqrt{K})$. \ Then for $x_{0}\in\Phi^{-1}(0)$
\[
\left\vert B_{k}r\right\vert ^{2}([x_{0}])=2^{-d/2}\mathrm{\operatorname*{vol}%
}(G\cdot x_{0})\,\left\vert r\right\vert ^{2}(x_{0}).
\]

\end{theorem}

The factor of $2^{-d/2}\mathrm{\operatorname*{vol}}(G\cdot x_{0})$ appearing
in Theorem \ref{thm:BKS comparison} is the ultimate reason that the modified
Guillemin--Sternberg map is asymptotically unitary; as we will see in\ Section
\ref{sec:asymptotics}, it precisely cancels the leading order asymptotic value
of the uncorrected density $I_{k}.$

\begin{proof}[Proof of Theorem \protect\ref{thm:modGS invertible}]
The natural map is injective because two holomorphic sections which agree on
the stable set $M_{s}$, which is an open dense subset of $M$, must
necessarily be equal. \ We now establish surjectivity for large $k.$
A section $\widehat{r}\in\mathcal{H}(M/\!\!/G;\hat{\ell}^{\otimes k}\otimes%
\sqrt{\widehat{K}})$ determines a $G$-invariant section $\pi^{\ast }\widehat{%
r}$ over the zero-set $\Phi^{-1}(0).$ Using the action of $G_{\mathbb{C}}$
on the bundles, we can extend $\pi^{\ast}\widehat{r}$ uniquely to a
holomorphic, $G_{\mathbb{C}}$-invariant section of $\ell^{\otimes k}\otimes%
\sqrt{K}$ defined over the stable set $M_{s}.$ As shown in
\cite[Appendix]{Guillemin-Sternberg}, there exists some $k$ for which $\mathcal{H}%
(M;\ell^{\otimes k})$ contains some nonzero $G$-invariant section $s.$ By %
(\ref{eqn:stable set character}), the unstable set is contained in the zero-set of $s.$ Thus, if the magnitude of $r$ remains bounded as we approach the
unstable set, the Riemann Extension Theorem (e.g.,\cite[pp. 9]{Griffiths})
will imply that $r$ extends holomorphically to all of $M.$
Suppose $\xi\in\mathfrak{g}$ has $\left\vert \xi\right\vert =1.$ Then we
will show in the next section (Theorem \ref{thm:s along a curve}) that the
variation of the magnitude $\left\vert r\right\vert ^{2}$ along the curve $%
e^{it\xi}\cdot x_{0}$ is
\begin{equation*}
\frac{d}{dt}\left\vert r\right\vert ^{2}(e^{it\xi}\cdot x_{0})=\left\vert
r\right\vert ^{2}(x_{0})\left[ -2k\,\phi_{\xi}(e^{it\xi}\cdot x_{0})-\frac{%
\mathfrak{L}_{JX^{\xi}}\varepsilon_{\omega}}{2\varepsilon_{\omega}}%
(e^{it\xi}\cdot x_{0})\right]
\end{equation*}
where $\varepsilon_{\omega}=\omega^{\wedge n}/n!$ is the Liouville volume
form. Furthermore, we will show in (\ref{gradient}) that $%
\phi_{\xi}(e^{it\xi }\cdot x_{0})$ increases with $t$ for $t\geq0.$ As a
result, the function $\phi_{\xi}(e^{i\xi}\cdot x_{0}),$ as $\xi$ varies over
unit vectors in $\mathfrak{g}$ and $x_{0}$ varies over $\Phi^{-1}(0),$ is
strictly positive and thus bounded below by compactness. Meanwhile, $%
\mathfrak{L}_{JX^{\xi}}\varepsilon_{\omega}/\varepsilon_{\omega}$ is bounded
over $M$ uniformly in $\xi$ with $\left\vert \xi\right\vert =1,$ again by
compactness. Using the monotonicity of $\phi_{\xi}(e^{it\xi}\cdot x_{0})$ in
$t,$ it follows that for all sufficiently large $k,$%
\begin{equation*}
2k\phi_{\xi}(e^{it\xi}\cdot x_{0})\geq\left\vert \frac{\mathfrak{L}%
_{JX^{\xi}}\varepsilon_{\omega}}{2\varepsilon_{\omega}}(e^{it\xi}\cdot
x_{0})\right\vert
\end{equation*}
for all $x_{0}\in\Phi^{-1}(0)$, all $\xi\in\mathfrak{g}$ with $\left\vert
\xi\right\vert =1,$ and all $t\geq1.$ It follows that if $k$ is large
enough, every $r$ obtained in the above way will extend holomorphically to
all of $M.$
\end{proof}

The proof of Theorem \ref{thm:BKS comparison} boils down to the following two
lemmas. Recall that $\Xi=\{\xi_{1},\dots,\xi_{d}\}$ is a basis of
$\mathfrak{g}$ for which the associated Haar measure on $G$ is normalized to
$1$.

\begin{lemma}
\label{lemma:det = vol} The function $\sqrt{\det_{\Xi}B_{x}}:=\det
(B(X^{\xi_{j}},X^{\xi_{k}}))$ is constant along each $G$-orbit. Moreover,
\[
\sqrt{\det{}_{\Xi}B_{x}}=\mathrm{\operatorname*{vol}}(G\cdot x).
\]

\end{lemma}

\begin{proof}
The fixed basis $\Xi$ of ${\mathfrak{g}}$ defines a left-invariant coframe $%
\vartheta$ on $G\cdot x.$ With respect to this coframe, the Riemannian
volume associated to $B$ is
\begin{equation*}
d\mathrm{\operatorname{vol}}=\sqrt{\det{}_{\Xi}B}\,\vartheta^{1}\wedge
\dots\wedge\vartheta^{k}.
\end{equation*}
By our definition of $\Xi,$ the pullback $dg_{0}=\varphi^{\ast}\vartheta
^{1}\wedge\dots\wedge\vartheta^{k}$ to $G$ by the map $\varphi:G\rightarrow
G\cdot x$ is a left Haar measure on $G$ for which the corresponding volume
is $\mathrm{\operatorname{vol}}_{0}(G)=1.$ Denote by $d\mathrm{\operatorname{vol}}%
_{0}(G)=\vartheta^{1}\wedge\dots\wedge \vartheta^{k}.$ Then since our choice
of inner product on $\mathfrak{g}$ yields $\int_{G\cdot x}d\mathrm{\operatorname{%
vol}}_{0}(G)=1$, integrating the equation $d\mathrm{\operatorname{vol}}=\sqrt{%
\det_{\Xi}B}\ d\mathrm{\operatorname{vol}}_{0}(G) $ over $G\cdot x$ yields the
desired result.
\end{proof}

Observe that if we choose a basis $\Xi$ which is not orthonormal with respect
to our fixed normalized inner product on $G$, or if we choose an inner product
on $\mathfrak{g}$ that is not normalized so that $\int_{G}dg=1$, then Lemma
\ref{lemma:det = vol} yields $\sqrt{\det_{\Xi}B_{x}}%
=C\mathrm{\operatorname*{vol}}(G\cdot x).$ The final effect of our seemingly
natural choices (which make $C$ equal $1$) is that the density $J_{k}$ (which
takes into account the metaplectic correction) approaches $1$ as
$k\rightarrow\infty$ (in general it would be $C^{-1}$).

Our second lemma is technical; it is a straightforward though somewhat tedious calculation.

\begin{lemma}
Let $Z^{j}=\pi_{+}X^{\xi_{j}}$. Then for $x_{0}\in\Phi^{-1}(0)$ we have%
\[
\left.  \mathfrak{i}(\bigwedge\nolimits_{j}Z^{j})\circ\mathfrak{i}%
(\bigwedge\nolimits_{k}\bar{Z}^{j})\,\varepsilon_{\omega}(x_{0})\right\vert
_{\Phi^{-1}(0)}=2^{-d}\left(  \mathrm{\operatorname*{vol}}(G\cdot
x_{0})\right)  ^{2}\left.  \,\frac{\omega^{n-d}}{(n-d)!}(x_{0})\right\vert
_{\Phi^{-1}(0)}.
\]

\end{lemma}

\begin{proof}
First, since $M$ is K\"ahler, the symplectic form $\omega$ is of $(1,1)$%
-type, that is, $\omega$ contracted on two holomorphic or two
antiholomorphic vectors is zero.
With our convention that $\bigwedge\nolimits_{j}w^{j}:=\frac{1}{d!}%
\sum\nolimits_{\sigma \in S_{d}}w^{\sigma (1)}\otimes \cdots \otimes
w^{\sigma (d)}$ we get $\mathfrak{i}(\bigwedge_{j}w^{j})\alpha =\alpha
(w^{1},\dots ,w^{d},\cdot ,\dots ,\cdot )=\mathfrak{i}(w^{d})\circ \cdots
\circ \mathfrak{i}(w^{1})\alpha .$ Using the fact that the interior product
is an antiderivation, we compute%
\begin{equation*}
\mathfrak{i}(\bar{Z}^{d})\circ \cdots \circ \mathfrak{i}(\bar{Z}^{1})\omega
^{n}=(-1)^{d(d-1)/2+d}\frac{n!}{(n-d)!}\bigwedge\nolimits_{j}\omega (\cdot ,%
\bar{Z}^{j})\wedge \omega ^{n-d}.
\end{equation*}
Continuing on, we contract the above with the holomorphic polyvector:
\begin{align}
\mathfrak{i}(\bigwedge\nolimits_{j}Z^{j}{)\circ} {\mathfrak{i}(}\bigwedge
\nolimits_{k}\bar{Z}^{k}\,){\text{\thinspace}}\,\omega^{n} &=(-1)^{d(d+1)/2}%
\frac{n!}{(n-d)!}\big[\det(\omega(Z^{j},\bar{Z}^{k}))\omega^{n-d}  \notag \\
&\qquad+(-1)^{d}\bigwedge\nolimits_{j}\omega(\cdot,\bar{Z}%
^{j})\wedge(n-d)!\bigwedge\nolimits_{k}\,\omega(Z^{k},\cdot)\wedge\omega
^{n-2d}\big].  \label{eqn:tech lemma 1}
\end{align}
A short computation shows that $\omega(Z^{j},\bar{Z}^{k})=\frac{i}{2}%
B(X^{\xi_{j}},X^{\xi_{k}}).$ When restricted to $\Phi^{-1}(0),$ the moment
map is constant. Hence
\begin{equation*}
\left. \omega(\cdot,\pi_{\pm}X^{\xi_{j}})\right\vert _{\Phi^{-1}(0)}=\left.
\tfrac{1}{2}d\phi_{\xi_{j}}\mp\tfrac{i}{2}\omega(\cdot,JX^{\xi_{j}})\right\vert _{\Phi^{-1}(0)}=\left. \mp\tfrac{i}{2}\omega(\cdot,JX^{\xi_{j}})\right\vert _{\Phi^{-1}(0)}
\end{equation*}
and so $\left. \bigwedge\nolimits_{j}\omega(\cdot,\bar{Z}^{j})\wedge
\bigwedge\nolimits_{k}\omega(Z^{k},\cdot)\right\vert _{\Phi^{-1}(0)}=0.$
Equation (\ref{eqn:tech lemma 1}) is now
\begin{equation*}
\left. \mathfrak{i}(\bigwedge\nolimits_{j}Z^{j})\circ\mathfrak{i}%
(\bigwedge\nolimits_{k}\bar{Z}^{k}\,){\text{\thinspace}}\,\frac{\omega^{n}}{%
n!}\right\vert _{\Phi^{-1}(0)}=\left. 2^{-d}(\det{}_{\Xi}B)\,\frac {%
\omega^{n-d}}{(n-d)!}\,\right\vert _{\Phi^{-1}(0)}.
\end{equation*}
Combining this with Lemma \ref{lemma:det = vol} yields the desired result.
\end{proof}

\begin{proof}[Proof of Theorem \protect\ref{thm:BKS comparison}]
Near $x_{0}$ we can write $r=s\nu$ where $s$ is a local $G$-invariant
holomorphic section of $\ell^{\otimes k}$ and $\nu$ is a local $G$-invariant
holomorphic section of $\sqrt{K}$. Let $\alpha=\nu^{2}.$ Then since $%
\mathfrak{i}(\pi_{+}X^{\xi})\,\bar{\alpha }=\mathfrak{i}(\pi_{-}X^{\xi})\,%
\alpha=0,$ we have
\begin{align*}
\pi^{\ast}\left( \left\vert B\nu\right\vert ^{2}\varepsilon_{\hat{\omega}%
}\right) & =\left. \left( \pi^{\ast}B\nu,\pi^{\ast}B\nu\right) ^{2}\,\frac{%
\omega^{n-d}}{(n-d)!}\right\vert _{\Phi^{-1}(0)} \\
&=\,\left. \mathfrak{i}(\bigwedge\nolimits_{j}Z^{j}\,\text{{\thinspace}}%
)\,\alpha \wedge\mathfrak{i}(\bigwedge\nolimits_{k}\bar{Z}^{k}{)}\,\bar{%
\alpha }\right\vert _{\Phi^{-1}(0)} \\
&=\left. \,\mathfrak{i}(\bigwedge\nolimits_{j}Z^{j}\,{\wedge}%
\bigwedge\nolimits_{k}\bar{Z}^{k}\,{)}((\nu,\nu)^{2}\frac{\omega^{n}}{n!}%
)\right\vert _{\Phi^{-1}(0)} \\
&=\,\left. \left( \nu,\nu\right) ^{2}~2^{-d}\det{}_{\Xi}B\ \frac{\omega^{n-d}%
}{(n-d)!}\right\vert _{\Phi^{-1}(0)}
\end{align*}
The result now follows upon dividing by $\omega^{n-d}/(n-d)!$, taking the
square root and using Lemma \ref{lemma:det = vol} and the fact that $%
\left\vert A_{k}s\right\vert ^{2}([x_{0}])=\left\vert s\right\vert
^{2}(x_{0})$.
\end{proof}

\section{Norm decompositions\label{sec:fiber integration}}

In this section, we show how the norm-squared of a $G$-invariant holomorphic
section over the \textit{upstairs} manifold $M$ can be computed as an integral
over the \textit{downstairs} manifold $M/\!\!/G.$ To do this, we first show
that the pointwise magnitude of a $G$-invariant holomorphic section (with or
without the metaplectic correction) varies in a predictable way as we move off
the zero-set of the moment map along the orbits of $G_{\mathbb{C}}.$ This
means that the behavior of an invariant section on $\Phi^{-1}(0)$ determines
its behavior on the whole stable set $M_{s}.$  We then integrate the resulting
expressions over the stable set using the decomposition in Theorem \ref{cor:lambda}.
By (\ref{eqn:stable set character}), the unstable set is contained in a set of
complex codimension at least $1$ and hence the stable set is a set of full
measure in $M.$ A similar computation shows that a Toeplitz operator on the upstairs
manifold with a $G $-invariant symbol, when restricted to the space of $G$-invariant
sections, is equivalent to a certain Toeplitz on the downstairs manifold.

Recall from Section \ref{subsec:Kahler reduction} that the stable set $M_{s}$
consists of those points in $M$ that can be moved into the zero-set of the
moment map by means of the action of $G_{\mathbb{C}}.$ Recall also (Theorem
\ref{cor:lambda}) that every point in $M_{s}$ can be expressed uniquely in the
form $e^{i\xi}\cdot x_{0},$ for some $\xi\in\mathfrak{g}$ and $x_{0}\in
\Phi^{-1}(0).$ The first main result of this section shows how the (pointwise)
magnitude of a $G$-invariant holomorphic section, evaluated at $e^{i\xi}\cdot
x_{0}$, varies with respect to $\xi.$

\begin{theorem}
\label{thm:s along a curve}Let $s$ be a $G$-invariant holomorphic section of
$\ell^{\otimes k}$ and let $r$ be a $G$-invariant holomorphic section of
$\ell^{\otimes k}\otimes\sqrt{K}$. Let $x_{0}$ be a point in $\Phi^{-1}(0)$
and $\xi$ an element of $\mathfrak{g}.$ Then the magnitudes of the sections at
$e^{i\xi}\cdot x_{0}$ are related to the magnitudes at $x_{0}$ as follows:%
\begin{align}
\text{(a)} &  \left\vert s\right\vert ^{2}(e^{i\,\xi}\cdot x_{0})=\left\vert
s\right\vert ^{2}(x_{0})\exp\left\{  -\int_{0}^{1}2k\phi_{\xi}(e^{it\xi}\cdot
x_{0})~dt\right\} \label{eqn:s norm}\\
\text{(b)} &  \left\vert r\right\vert ^{2}(e^{i\xi}\cdot x_{0})=\left\vert
r\right\vert ^{2}(x_{0})\exp\Bigg\{-\int_{0}^{1}\Bigg(2k\phi_{\xi}(e^{it\xi
}\cdot x_{0})+\frac{{\mathcal{L}}_{JX^{\xi}}\varepsilon_{\omega}}%
{2\varepsilon_{\omega}}(e^{it\xi}\cdot x_{0})\Bigg)dt\Bigg\}\label{eqn:r norm}%
\end{align}
Furthermore, if we consider
\[
\rho(\xi,x_{0}):=2\int_{0}^{1}\phi_{\xi}(e^{it\xi}\cdot x_{0})~dt,
\]
then for each fixed $x_{0},$ $\rho(\xi,x_{0})$ achieves its unique minimum at
$\xi=0$ and the Hessian of $\rho(\xi,x_{0})$ at $\xi=0$ is given by%
\[
\left.  D_{\xi_{1}}D_{\xi_{2}}\rho(\xi,x_{0})\right\vert _{\xi=0}=2B_{x_{0}%
}(JX^{\xi_{1}},JX^{\xi_{2}}),\quad\xi_{1},\xi_{2}\in\mathfrak{g.}%
\]

\end{theorem}

Note that the extra factor in (\ref{eqn:r norm}), as compared to
(\ref{eqn:s norm}), is independent of $k$ and is equal to $1$ when $\xi$ is
equal to $0$. This extra factor, therefore, does not affect the leading order
asymptotics. On the other hand, this extra factor can become unbounded near
the unstable set which means that pointwise behavior of $r$ tends to be worse
than that of $s$.

Let $\Lambda:\mathfrak{g}\times\Phi^{-1}(0)\rightarrow M_{s}$ denote the
diffeomorphism $\Lambda(\xi,x_{0})=e^{i\xi}\cdot x_{0}$ of Theorem
\ref{cor:lambda}. Recall also that we have chosen an orthonormal basis
$\Xi=\{\xi_{j}\}_{j=1}^{d}$ of $\mathfrak{g}$ to which there corresponds a
Lebesgue measure $d^{d}\xi$ on $\mathfrak{g}$. The Liouville volume
$\varepsilon_{\omega}=\omega^{\wedge n}/n!$ on $M$ (which is the same as the
Riemannian volume) decomposes as%
\[
\Lambda^{\ast}(\varepsilon_{\omega})_{(\xi,x_{0})}=\tau(\xi,x_{0})\,d^{d}%
\xi\wedge\,d\mathrm{\operatorname*{vol}}(\Phi^{-1}(0))_{x_{0}}%
\]
for some $G$-invariant smooth Jacobian function $\tau\in C^{\infty
}(\mathfrak{g}\times\Phi^{-1}(0))$. We are now ready to state the remaining
main results of this section. Let $\gamma_{\xi}:[0,1]\rightarrow M_{s}$ be the
path $\gamma_{\xi}(t)=e^{it\xi}\cdot x_{0}$; then $\int_{\gamma_{\xi}}%
\phi_{\xi}=\int_{0}^{1}\phi_{\xi}(e^{it\xi}\cdot x_{0})~dt.$

\begin{theorem}
\label{thm:s norm decomp}Suppose $s$ is a $G$-invariant holomorphic section of
$\ell^{\otimes k}$. \ Then the norm of $s$ can be computed by%
\begin{equation}
\left\Vert s\right\Vert ^{2}:=\left(  k/2\pi\right)  ^{n/2}\int_{M}%
|s|^{2}\,\,\varepsilon_{\omega}=\left(  k/2\pi\right)  ^{(n-d)/2}%
\int_{M/\!\!/G}|A_{k}s|^{2}\,I_{k}\,~\varepsilon_{\widehat{\omega}%
}\label{first}%
\end{equation}
where, with $[x_{0}]=G\cdot x_{0}$,
\begin{equation}
I_{k}([x_{0}])=\mathrm{\operatorname*{vol}}(G\cdot x_{0})\,\left(
k/2\pi\right)  ^{d/2}\int_{\mathfrak{g}}\tau(\xi,x_{0})\,\exp\left\{
-2k\int_{\gamma_{\xi}}\phi_{\xi}\right\}  \,d^{d}\xi.\label{eqn:Ik}%
\end{equation}

\end{theorem}

In the case of a Hamiltonian torus action (i.e., the case when $G$ is
commutative), a similar formula was obtained by Charles in \cite[Sec
4.5]{Charles}. Nevertheless, our result is slightly stronger than that of
Charles, in that Charles inserts into the first integral in (\ref{first}) a
function $f$ that is assumed to be supported away from the unstable set. In
general, this allows him to obtain asymptotics about Toeplitz operators (as we
do also in Theorem \ref{toeplitz.thm}). To obtain results about the norm of a
section, Charles takes $f$ to be identically equal to one in a neighborhood of
the zero-set, but still supported away from the unstable set. The insertion of
such a cutoff function leaves unchanged the asymptotics to all orders of the
first integral in (\ref{first}). Still, our formula actually gives an exact
(not just asymptotic) expression for the norm of an invariant section upstairs
as an integral over the downstairs manifold. To get this, we require some
estimates for the blow-up of certain Jacobians as we approach the unstable set.

\begin{theorem}
\label{thm:half-form norm decomp}Suppose $r$ is a $G$-invariant holomorphic
section of $\ell^{\otimes k}\otimes\sqrt{K}.$ \ Then the norm of $r$ can be
computed by%
\[
\Vert r\Vert^{2}:=\left(  k/2\pi\right)  ^{n/2}\int_{M}|r|^{2}\,\varepsilon
_{\omega}=\left(  k/2\pi\right)  ^{(n-d)/2}\int_{M/\!\!/G}|B_{k}r|^{2}%
\,J_{k}\,\,\varepsilon_{\widehat{\omega}}%
\]
where
\begin{equation}
J_{k}([x_{0}])=\left(  k/2\pi\right)  ^{d/2}2^{d/2}\int_{\mathfrak{g}}\tau
(\xi,x_{0})\,\exp\left\{  -\int_{\gamma_{\xi}}\left(  2k\,\phi_{\xi}%
+\frac{{\mathcal{L}}_{JX^{\xi}}\varepsilon_{\omega}}{2\varepsilon_{\omega}%
}\right)  \right\}  \,d^{d}\xi.\label{eqn:Jk}%
\end{equation}

\end{theorem}

Note the volume factor that is present in the expression for $I_{k}$ but not
in the expression for $J_{k}.$ In computing $I_{k},$ we decompose $M_{s}$ as
$\mathfrak{g}\times\Phi^{-1}(0)$ and then integrate out the $\xi$-dependent
part of (\ref{eqn:s norm}), leaving an integration over $\Phi^{-1}(0).$ Since
the resulting integrand on $\Phi^{-1}(0)$ will be $G$-invariant, the
integration over $\Phi^{-1}(0)$ can be turned into an integration over
$\Phi^{-1}(0)/G.$ The volume factor in (\ref{eqn:Ik}) arises because the
projection map $\pi:\Phi^{-1}(0)\rightarrow\Phi^{-1}(0)/G$ is a Riemannian
submersion, whence the Riemannian volume measure on $\Phi^{-1}(0)$ maps to the
Riemannian volume measure on $\Phi^{-1}(0)/G, $ \textit{multiplied by} a
density given by the volume factor. This crucial fact is the basic reason that
the Guillemin--Sternberg map (without half-forms) is \emph{not} asymptotically unitary.

Meanwhile, the volume factor \textit{fails} to arise in (\ref{eqn:Jk}) because
it is canceled by the volume factor in Theorem \ref{thm:BKS comparison}, which
relates the pointwise magnitude of $r$ at $x_{0}$ to the pointwise magnitude
of $B_{k}r$ at $[x_{0}].$

In Section \ref{sec:asymptotics}, we will calculate the asymptotic behavior of
the expressions for densities $I_{k}$ and $J_{k}$ as $k$ tends to infinity.
The calculation is done by means of Laplace's method and uses the Hessian
computed in Theorem \ref{thm:s along a curve}. Because the extra factor in the
expression for $\left\vert r\right\vert ^{2}(e^{i\xi}\cdot x_{0})$ is
independent of $k$ and equal to 1 at $\xi=0,$ this factor does not affect the
leading-order asymptotics of the norm. The different asymptotic behavior of
the maps $A_{k}$ and $B_{k}$ is due to the volume factor in the expression for
$I_{k}$ that is not present in the expression for $J_{k}.$

If it should happen that $I_{k}$ is constant, then $A_{k}$ (which is in any
case invertible for all $k$) will be a constant multiple of a unitary map.
Similarly, if $J_{k}$ is constant for some large $k$ (large enough that
$B_{k}$ is invertible), then $B_{k}$ will be a constant multiple of a unitary
map. There does not, however, seem to be any reason that either $I_{k}$ or
$J_{k}$ should typically be constant. Nevertheless, we will see in Section
\ref{sec:asymptotics} that $J_{k}$ asymptotically approaches 1 for large $k,$
which implies that $B_{k}$ is asymptotically unitary. On the other hand,
$I_{k}$ is not asymptotically constant unless all the $G$-orbits in $\Phi
^{-1}(0)$ happen to have the same volume. In the case where the $G$-orbits do
not all have the same volume, it is not hard to show (Theorem
\ref{thm:non-unitarity}) that $A_{k}$ is not asymptotic to any constant
multiple of a unitary map.

A similar analysis shows that if we consider a Toeplitz operator with a
$G$-invariant symbol $f$ upstairs, then the matrix entries for such an
operator can be expressed an integral over the downstairs manifold involving a
certain density $I_{k}(f),$ which reduces to $I_{k}$ when $f\equiv1.$ Theorem
\ref{toeplitz.thm} shows that in the case with half forms we obtain an
asymptotic equivalence of a very simple form between Toeplitz operators
upstairs and downstairs.

\bigskip

We now turn to the proofs of Theorems \ref{thm:s along a curve} and
\ref{thm:s norm decomp}. First, we will decompose the integral over $M_{s}$
into integrals over $\mathfrak{g}$ and $\Phi^{-1}(0)$ and the integral over
$\Phi^{-1}(0)$ into integrals over $M/\!\!/G$ and $G\cdot x_{0}.$ Both of
these integral decompositions follow from the coarea formula \cite[pp. 159--160]{Chavel}:

\begin{lemma}
[Coarea formula]\label{lemma:coarea} Let $Q$ and $N$ be smooth Riemannian
manifolds with $\dim Q\geq\dim N$, and let $p:Q\rightarrow N$. Then for any
$f\in L^{1}(M)$ one has
\[
\int_{Q}\mathcal{J}_{p}\,f\,\,d\mathrm{\operatorname*{vol}}(Q)=\int
_{N}d\mathrm{\operatorname*{vol}}(N)(y)\int_{p^{-1}(y)}(\left.  f\right\vert
_{p^{-1}(y)})\,d\mathrm{\operatorname*{vol}}(p^{-1}(y))
\]
where the Jacobian is $\mathcal{J}_{p}:=\sqrt{\det p_{\ast}\circ p_{\ast
}^{\operatorname*{adj}}}.$
\end{lemma}

Applying the coarea formula to the map $\operatorname*{pr}\nolimits_{2}%
\circ\Lambda^{-1}:e^{i\xi}\cdot x_{0}\mapsto x_{0}$ from $M_{s}$ to $M_{0}$,
and identifying $\exp(i\mathfrak{g)}\cdot x_{0}$ with $\mathfrak{g}$, we have
that for every $f\in L^{1}(M_{s})$%
\[
\int_{M_{s}}f\,d\mathrm{\operatorname*{vol}}(M_{s})=\int_{M_{0}}%
\int_{\mathfrak{g}}(\left.  f\right\vert _{\exp(i\mathfrak{g)}\cdot x_{0}%
})\,\tau\,d^{d}\xi\,d\mathrm{\operatorname*{vol}}(M_{0})
\]
where $\tau$ is the Jacobian $\mathcal{J}$ of the map $\Lambda:\mathfrak{g}%
\times\Phi^{-1}(0)\rightarrow M_{s}.$

Since the volume form on the quotient $M/\!\!/G$ is $d$%
\textrm{$\operatorname*{vol} $}$(M/\!\!/G)=\widehat{\omega}\,^{n-d}/(n-d)!$
and since $\Phi^{-1}(0)\rightarrow M/\!\!/G$ is a Riemannian submersion, we
have \cite[Sec. 2.2]{Gilkey}%
\[
d\mathrm{\operatorname*{vol}}(\Phi^{-1}(0))=d\mathrm{\operatorname*{vol}%
}(G\cdot x_{0})\wedge\frac{\pi_{hor}^{\ast}\widehat{\omega}\,^{n-d}}{(n-d)!}%
\]
where $\pi_{hor}^{\ast}\widehat{\omega}\,^{n-d}$ denotes the horizontal lift
(i.e.,pullback composed with projection to the horizontal subspace $W_{x_{0}%
}^{\bot_{B}}$) of the Liouville form on the quotient. \ The $2$-form
$\pi^{\ast}\widehat{\omega}=i^{\ast}\omega$ is already horizontal, though,
since%
\[
\mathfrak{i}(X^{\xi})\,\pi^{\ast}\widehat{\omega}=\mathfrak{i(}X^{\xi
})\,i^{\ast}\omega=d\phi_{\xi}|_{\Phi^{-1}(0)}=0.
\]
This proves the following.

\begin{lemma}
\label{lemma:zero-set integral decomp}Let $\pi:\Phi^{-1}(0)\rightarrow
M/\!\!/G:=\Phi^{-1}(0)/G$ denote the canonical projection. Then $\mathcal{J}%
_{\pi}=1;$ in particular, for every $G$-invariant function $f\in L^{1}%
(\Phi^{-1}(0))^{G}$ we have%
\[
\int_{\Phi^{-1}(0)}f\,d\mathrm{\operatorname*{vol}}(\Phi^{-1}(0))=\int
_{M/\!\!/G}\mathrm{\operatorname*{vol}}(G\cdot x_{0})\,f([x_{0}%
])\,d\mathrm{\operatorname*{vol}}(M/\!\!/G)\,.
\]

\end{lemma}

To prove Theorem \ref{thm:s along a curve}, we first make the following computations:

\begin{lemma}
\label{lemma:derivatives}The first two derivatives of the norm of a
$G$-invariant holomorphic section $s\in\mathcal{H}(M;\ell^{\otimes k})^{G}$
along $\gamma_{\xi}$ are:%
\begin{align*}
&  JX^{\xi}\left\vert s\right\vert ^{2}=-2k\,\phi_{\xi}\left\vert s\right\vert
^{2},\ {\text{and}}\\
&  JX^{\xi_{1}}JX^{\xi_{2}}\left\vert s\right\vert ^{2}=-2k~B(X^{\xi_{1}%
},X^{\xi_{2}})\left\vert s\right\vert ^{2}+4k^{2}~\phi_{\xi_{1}}\phi_{\xi_{2}%
}\left\vert s\right\vert ^{2}.
\end{align*}

\end{lemma}

The last statement of Theorem \ref{thm:s along a curve} follows immediately
upon combining the two above equations and restricting to the zero-set.

\begin{proof}
Since the connection on $\ell^{\otimes k}$ is Hermitian, we have
\begin{equation*}
JX^{\xi}\left\vert s\right\vert ^{2}=(\nabla_{JX^{\xi}}\,s,s)+(s,\nabla
_{JX^{\xi}}\,s).
\end{equation*}
By assumption, $s$ is $G$-invariant. Since the action of $G$ on $\mathcal{H}%
(M;\ell^{\otimes k})$ is given by equation (\ref{eqn:G action on sections}),
we have
\begin{equation*}
0=Q_{\xi}s=\left( \nabla_{X^{\xi}}-ik\,\,\phi_{\xi}\right) s
\end{equation*}
so that $\nabla_{X^{\xi}}\,s=-ik\,\phi_{\xi}\,s.$ The projection of a vector
field onto its $(0,1)$-part is $X\mapsto\tfrac{1}{2}(1-iJ)X.$ We assume that
$s$ is holomorphic, i.e.,$\nabla_{(1-iJ)X^{\xi}}\,s=0.$ Therefore $%
\nabla_{JX^{\xi}}s=i\nabla_{X^{\xi}}s.$ Putting this all together, we get
\begin{equation*}
JX^{\xi}(s,s)=-i\left( ik\,\phi_{\xi}s,s\right) +i\left( s,ik\,\phi_{\xi
}s\right) =-2k\,\phi_{\xi}(s,s).
\end{equation*}
Next, the Leibniz rule yields the second derivative:
\begin{equation*}
JX^{\xi_{j}}JX^{\xi_{l}}\left\vert s\right\vert
^{2}=-2k\,JX^{\xi_{j}}(\phi_{\xi_{l}})+4k^{2}\,\phi_{\xi_{l}}\phi_{\xi_{j}}%
\left\vert s\right\vert ^{2}.
\end{equation*}
Since $d\phi_{\xi_{l}}=\mathfrak{i}(X^{\xi_{l}})\,\omega$ we have
\begin{equation*}
JX^{\xi_{j}}\phi_{\xi_{l}}=\omega(X^{\xi_{l}},JX^{\xi_{j}})=B(X^{%
\xi_{j}},X^{\xi_{l}}).
\end{equation*}
Therefore,
\begin{equation*}
JX^{\xi_{j}}JX^{\xi_{l}}\left\vert s\right\vert
^{2}=-2k\,B(X^{\xi_{j}},X^{\xi_{l}})\left\vert s\right\vert
^{2}+4k^{2}\,\phi_{\xi_{j}}\phi_{\xi_{l}}\left\vert s\right\vert ^{2}
\end{equation*}
as desired.
\end{proof}

\begin{proof}[Proof of Theorem \protect\ref{thm:s along a curve}(a)]
Applying Lemma \ref{lemma:derivatives}, the derivative of $\left\vert
s\right\vert ^{2}$ along $\gamma_{\xi}$ is
\begin{equation*}
\frac{d}{dt}\gamma_{\xi}^{\ast}\left\vert s\right\vert ^{2}(t)=JX^{\xi
}\left\vert s\right\vert ^{2}(e^{it\xi}\cdot x_{0})=-2k\,\phi_{\xi}(e^{it\xi
}\cdot x_{0})\,\left\vert s\right\vert ^{2}(e^{it\xi}\cdot x_{0}).
\end{equation*}
so that%
\begin{equation*}
\frac{d}{dt}\log\gamma_{\xi}^{\ast}\left\vert s\right\vert ^{2}(t)=-2k\,\phi
(e^{it\xi}\cdot x_{0}).
\end{equation*}
Integrating this equation from $t=0$ to $t=1$ yields the desired result.
\end{proof}

As Guillemin and Sternberg pointed out, an important consequence of Theorem
\ref{thm:s along a curve} is that a $G$-invariant holomorphic $s\in
\mathcal{H}(M;\ell^{\otimes k})$ attains its (unique) maximum value in the
orbit $G_{\mathbb{C}}\cdot x_{0}$ on the zero-set of the moment map
$(G_{\mathbb{C}}\cdot x_{0})\cap\Phi^{-1}(0)=G\cdot x_{0}.$ The reason is that
$JX^{\xi}$ is the gradient vector of the function $\phi_{\xi}$, since%
\begin{equation}
B(JX^{\xi},\cdot)=\omega(JX^{\xi},J\cdot)=\omega(X^{\xi},\cdot)=d\phi_{\xi
}.\label{gradient}%
\end{equation}
Hence,$\ \phi_{\xi}$ is increasing along $\gamma_{\xi}$ and so $\left\vert
s\right\vert ^{2}$ is decreasing away from $\Phi^{-1}(0).$

A useful and direct consequence of this fact is that each $G_{\mathbb{C}}%
$-orbit intersects the zero-set in exactly one $G$-orbit; since if not, then
there exist points $x,x^{\prime}\in\Phi^{-1}(0)$ which lie in distinct
$G$-orbits such that $x^{\prime}=e^{i\xi}\cdot x$ for some $\xi\in
\mathfrak{g}$. But this is impossible since moment map is strictly increasing
along the path $t\mapsto e^{it\xi}\cdot x$.

\bigskip

Next we express the norm of a $G$-invariant half-form corrected section $r\in$
$\mathcal{H}(M;\ell^{\otimes k}\otimes\sqrt{K})^{G}$ along $\gamma_{\xi}$ in
terms of its value at $x_{0}$.

\begin{proof}[Proof of Theorem \protect\ref{thm:s along a curve}(b)]
Locally, a section $r\in\mathcal{H}(M;\ell^{\otimes k}\otimes\sqrt{K})^{G} $
can be written as $r=s\nu$ for some $G$-invariant holomorphic half-form $\nu$%
. Since $\nu$ is a $G$-invariant holomorphic half-form, the combination $%
\nu^{2}\wedge \overline{\nu}^{2}$ is $G_{\mathbb{C}}$-invariant, whence
\begin{equation*}
0={\mathcal{L}}_{JX^{\xi}}(\nu^{2}\wedge\overline{\nu}^{2})={\mathcal{L}}%
_{JX^{\xi}}(\left( \nu,\nu\right) ^{2}\varepsilon_{\omega}) =\left[ 2\left(
\nu,\nu\right) JX^{\xi}\left( \nu,\nu\right) \right] \,\varepsilon_{\omega}+%
\left( \nu,\nu\right) ^{2}{\mathcal{L}}_{JX^{\xi}}\varepsilon_{\omega}
\end{equation*}
so that
\begin{equation*}
\frac{d}{dt}\log\gamma_{\xi}^{\ast}\left( \nu,\nu\right) (t)=\frac{JX^{\xi
}\left( \nu,\nu\right) }{\left( \nu,\nu\right) }(e^{it\xi}\cdot x_{0})=-%
\frac{{\mathcal{L}}_{JX^{\xi}}\varepsilon_{\omega}}{2\varepsilon _{\omega}}%
(e^{it\xi}\cdot x_{0}).
\end{equation*}
Integrating this along $\gamma_{\xi}$ and combining it with Theorem \ref%
{thm:s along a curve}(a) yields the desired result.
\end{proof}

By decomposing the Liouville measure on $M_{s}$ in terms of the global
decomposition $\psi:\mathfrak{g}\times\Phi^{-1}(0)\rightarrow M_{s}$ (Theorem
\ref{cor:lambda}) and the fibration $G\rightarrow\Phi^{-1}(0)\rightarrow
M/\!\!/G$, we will obtain our desired integral. Recall that we have chosen a
basis $\Xi=\{\xi_{j}\}_{j=1}^{d}$ of $\mathfrak{g}$ to which there corresponds
a Lebesgue measure $d^{d}\xi$ on $\mathfrak{g}$.

The Liouville volume $\varepsilon_{\omega}=\omega^{\wedge n}/n!$, which is the
same as the Riemannian volume, decomposes as%
\begin{equation}
\Lambda^{\ast}(\varepsilon_{\omega})_{(\xi,x_{0})}=\tau(\xi,x_{0})\,d^{d}%
\xi\wedge\,d\mathrm{\operatorname*{vol}}(\Phi^{-1}(0))_{x_{0}}%
\label{eqn:liouville decomp}%
\end{equation}
for some $G$-invariant smooth Jacobian function $\tau\in C^{\infty
}(\mathfrak{g}\times\Phi^{-1}(0))$ (this is just the coarea formula
(\ref{lemma:coarea}) applied to the map $M_{s}\rightarrow M_{0}$ given by
$e^{i\xi}\cdot x_{0}\mapsto x_{0}$).

\begin{proof}[Proof of Theorem \protect\ref{thm:s norm decomp}]
Combining (\ref{eqn:liouville decomp}) and Theorem \ref{thm:s along a
curve} yields
\begin{multline}
\left( k/2\pi\right) ^{n/2}\int_{M}|s|^{2}\,\varepsilon_{\omega }
\label{eqn:f integral decomp} \\
=\left( k/2\pi\right) ^{(n-d)/2}\int_{M_{0}}\left( k/2\pi\right) ^{d/2}\int_{%
\mathfrak{g}}|\pi^{\ast}A_{k}s|^{2}(x_{0})\,\exp\left\{
-2k\int_{\gamma_{\xi}}\phi_{\xi}\right\} \\
\times\tau(\xi,x_{0})\ \,\,d^{d}\xi\,d\mathrm{\operatorname{vol}}(M_{0})
\end{multline}
since $s(x_{0})=\pi^{\ast}A_{k}s(x_{0})$.
The terms $|\pi ^{\ast }A_{k}s|^{2}$ and $\tau $ are obviously $G$%
-invariant. To see that the exponential factor in the integrand is $G$%
-invariant, define
\begin{equation*}
\rho (e^{i\xi }\cdot x_{0})=2\int_{\gamma _{\xi }}\phi _{\xi
}=2\int_{0}^{1}\left\langle \Phi (e^{it\xi }\cdot x_{0}),\xi \right\rangle
dt.
\end{equation*}%
To evaluate $g^{\ast }\rho $, we rewrite $ge^{i\xi }\cdot x_{0}=e^{iAd(g)\xi
}g\cdot x_{0}.$ Then%
\begin{align*}
g^{\ast }\rho (e^{i\xi }\cdot x_{0})& =\rho (e^{iAd(g)\xi }g\cdot
x_{0})=2\int_{0}^{1}\left\langle \Phi (e^{iAd(g)\xi t}g\cdot x_{0}),\operatorname{%
Ad}(g)\xi \right\rangle \,dt \\
& =2\int_{0}^{1}\left\langle \Phi (e^{it\xi }\cdot x_{0}),\xi \right\rangle
dt=\rho (e^{i\xi }\cdot x_{0})
\end{align*}%
where we have used the $G$-equivariance of the moment map $\Phi (gp)=%
\operatorname{Ad}^{\ast }(g^{-1})\Phi (p).$ Hence, the integrand in (\ref{eqn:f
integral decomp}) is $G$-invariant, so by the coarea formula (Lemma \ref%
{lemma:zero-set integral decomp}) the integral over $M_{0}$ reduces to an
integral over $M/\!\!/G$ with an extra factor of $\mathrm{\operatorname{vol}}%
(G\cdot x_{0}),$ and we obtain the desired result.
\end{proof}

We now compute the norm of a $G$-invariant holomorphic section as an integral
over the zero-set of the moment map taking into account the metaplectic correction.

\begin{proof}[Proof of Theorem \protect\ref{thm:half-form norm decomp}]
The proof is the similar to the above proof of Theorem \ref{thm:s norm
decomp}; we first write the integral over $M$ as a multiple integral:
\begin{multline*}
\left( k/2\pi\right) ^{n/2}\int_{M}|r|^{2}\,\varepsilon_{\omega}=\left(
k/2\pi\right) ^{(n-d)/2}\int_{M_{0}}\left( k/2\pi\right) ^{d/2} \\
\times\int_{\mathfrak{g}}|r|^{2}(x_{0})\exp\left\{
-\int_{\gamma_{\xi}}(2k\,\phi_{\xi}+{\mathcal{L}}_{JX^{\xi}}\varepsilon_{%
\omega}/2\varepsilon _{\omega})\right\} \  \\
\times\tau(\xi,x_{0})\ \,d^{d}\xi\,\,d\mathrm{\operatorname{vol}}(M_{0}).
\end{multline*}
The main difference from the proof of Theorem \ref{thm:s norm decomp} is
that by Theorem \ref{thm:BKS comparison} we have
\begin{equation*}
|r|^{2}(x_{0})=2^{d/2}\mathrm{\operatorname{vol}}(G\cdot
x_{0})^{-1}|B_{k}r|^{2}([x_{0}]).
\end{equation*}
Inserting this into the above integral, and noting that the remaining
integrand is $G$-invariant (following the same argument as in the proof of
Theorem \ref{thm:s norm decomp}) we obtain (again using the coarea formula
Lemma \ref{lemma:coarea}) that
\begin{multline*}
\left( k/2\pi\right) ^{n/2}\int_{M}|r|^{2}\,\varepsilon_{\omega}=\left(
k/2\pi\right) ^{(n-d)/2}\int_{M/\!\!/G}|r|^{2}([x_{0}])\,\  \\
\times\left( k/2\pi\right) ^{d/2}2^{d/2}\int_{\mathfrak{g}}\exp\left\{
-\int_{\gamma_{\xi}}(2k\,\phi_{\xi}+\frac{{\mathcal{L}}_{JX^{\xi}}%
\varepsilon_{\omega}}{2\varepsilon_{\omega}})\right\} \  \\
\times\tau(\xi,x_{0})\ d\mathrm{\operatorname{vol}}(\exp(i\mathfrak{g})\cdot
x_{0})\,d^{d}\xi\,\,\varepsilon_{\widehat{\omega}}
\end{multline*}
whence
\begin{equation*}
\Vert r\Vert^{2}=\left( k/2\pi\right)
^{(n-d)/2}\int_{M/\!\!/G}|B_{k}r|^{2}([x_{0}])J_{k}([x_{0}])\,\varepsilon_{%
\widehat{\omega}}
\end{equation*}
as desired.
\end{proof}

\section{Asymptotics\label{sec:asymptotics}}

In this section we compute the leading order asymptotics of the densities%
\[I_{k}(f)([x_{0}]) =\mathrm{\operatorname*{vol}}(G\cdot x_{0})(k/2\pi)^{d/2}
\int_{\mathfrak{g}}\tau(\xi,x_{0})\,f(x_{0},\xi)\,
\exp\left\{-2k\int_{\gamma_{\xi}}\phi_{\xi}\right\}  \,d^{d}\xi,\]
and
\[J_{k}(f)([x_{0}]) =(k/2\pi)^{d/2}\,2^{d/2}\int_{\mathfrak{g}}
\tau(\xi,x_{0})\,f(x_{0},\xi)\,\exp
\left\{-\int_{\gamma_{\xi}}\left(2k\,\phi_{\xi}+\frac{{\mathcal{L}}_{JX^{\xi}}\varepsilon_{\omega}}{2\varepsilon_{\omega}}\right)\right\}\,d^{d}\xi,\]
where $f$ is a smooth, $G$-invariant function on $M.$ In the case that
$f\equiv1$, these asymptotics will give asymptotics of the
Guillemin--Sternberg-type maps. For arbitrary smooth $f\in C^{\infty}(M)^{G}
$, the asymptotics will give information about Toeplitz operators.

The main result is the following.

\begin{theorem}
\label{thm:density asymptotics}For $f\in C^{\infty}(M)^{G},$ the densities
$I_{k}(f)$ and $J_{k}(f)$ satisfy
\begin{align*}
\lim_{k\rightarrow\infty}I_{k}(f)([x_{0}]) &  =2^{-d/2}f(x_{0}%
)\mathrm{\operatorname*{vol}}(G\cdot x_{0}),\text{ and }\\
\lim_{k\rightarrow\infty}J_{k}(f)([x_{0}]) &  =f(x_{0})
\end{align*}
for each $x_{0}\in\Phi^{-1}(0),$ and the limits are uniform.
\end{theorem}

Our main tool will be Laplace's approximation \cite[Chap. 10]{Bleistein}, also
frequently referred to as the stationary phase approximation or the method of
steepest descent. Let $D\subset\mathbb{R}^{d}$ be a bounded domain, and
consider $\rho\in C^{2}(D)$ and $\sigma\in C(D)$. Suppose $\rho$ attains a
unique minimum at $x_{0}\in D\setminus\partial D$ (i.e.,the interior of $D$).
\ Laplace's approximation gives the leading order asymptotic limit%
\begin{equation}
I(k)=\int_{D}\sigma(x)\,e^{-k\,\rho(x)}d^{d}x\sim e^{k\,\rho(x_{0})}\left(
\frac{2\pi}{k}\right)  ^{d/2}\left\vert \det\,H_{\rho}(x_{0})\right\vert
^{-1/2}\sigma(x_{0}),~k\rightarrow\infty,\label{eqn:Laplace's approx}%
\end{equation}
where\ $H_{\rho}$ denotes the Hessian of $\rho$.

The formula for the large-$k$ limits of $I_{k}$ and $J_{k}$ come from applying
Laplace's method to the integral over $\mathfrak{g}\simeq\mathbb{R}^{d}$ in
Theorems \ref{thm:s norm decomp} and \ref{thm:half-form norm decomp}, using the
computation of the relevant Hessian, which we have already performed in Theorem
\ref{thm:s along a curve}. Conceptually, then, it is easy to understand where
the limiting formulas come from.

There are, however, some technicalities to attend to, and it is these
technicalities that will occupy the bulk of this section. First, we need to
consider the part of the integrals near infinity as well as the part near the
origin. The density $\tau(x_{0},\xi)$ can (apparently) blow up as $\xi$ tends
to infinity for certain values of $x_{0}$. As a result, we need some estimates
(Lemma \ref{lemma:tau estimates}) to ensure that the $\xi$-integrals in
(\ref{eqn:Ik}) or (\ref{eqn:Jk}) are finite for all (as opposed to almost all)
$x_{0},$ at least for large $k.$ Second, we wish to be careful in verifying
that the limits are uniform over $M/\!\!/G,$ which is needed to obtain the
asymptotic unitarity of the maps $B_{k}.$

Before coming to these technicalities we state the consequences of the above
asymptotic formulas for the unitarity of the maps $A_{k}$ and $B_{k}$ and
indicate some applications to the asymptotics of Toeplitz operators.

\begin{theorem}
\label{thm:Bk unitary}The maps $B_{k}$ are asymptotically unitary, in the
sense that
\[
\lim_{k\rightarrow\infty}\left\Vert B_{k}^{\ast}B_{k}-I\right\Vert
=\lim_{k\rightarrow\infty}\left\Vert B_{k}B_{k}^{\ast}-I\right\Vert =0,
\]
where $\left\Vert \cdot\right\Vert $ refers to the operator norm.
\end{theorem}

\begin{proof}
Let $\widehat{r}\in\mathcal{H}(M/\!\!/G;\hat{\ell}^{\otimes k}\otimes\sqrt {\widehat{K}%
}).$ Then for all $t\in\mathcal{H}(M;\ell^{\otimes k}\otimes\sqrt{K})$ we
have
\begin{equation*}
\left\langle B_{k}^{\ast}\widehat{r},t\right\rangle _{M}=\left\langle
\widehat{r},B_{k}t\right\rangle _{M/\!\!/G}.
\end{equation*}
On the other hand, Theorem \ref{thm:half-form norm decomp} implies%
\begin{equation*}
\left\langle B_{k}^{\ast}\widehat{r},t\right\rangle _{M}=\left\langle
J_{k}B_{k}B_{k}^{\ast}\widehat{r},B_{k}t\right\rangle _{M/\!\!/G},~~\forall
t\in\mathcal{H}(M;\ell^{\otimes k}\otimes\sqrt{K}).
\end{equation*}
Since $B_{k}$ is bijective (for $k$ sufficiently large), every $\widehat{t}%
\in\mathcal{H}(M/\!\!/G;\hat {\ell}^{\otimes k}\otimes\sqrt{\widehat{K}})$ is of the
form $B_{k}t$ for some $t\in\mathcal{H}(M;\ell^{\otimes k}\otimes\sqrt{K});$
so for all $\widehat {t}\in\mathcal{H}(M/\!\!/G;\hat{\ell}^{\otimes k}\otimes%
\sqrt{\widehat{K}})$ we have%
\begin{equation*}
\left\langle \widehat{r},\widehat{t}\right\rangle _{M/\!\!/G}=\left\langle
J_{k}B_{k}B_{k}^{\ast}\widehat{r},\widehat{t}\right\rangle _{M/\!\!/G}.
\end{equation*}
Applying Theorem \ref{thm:density asymptotics} then yields the claim that $%
B_{k}$ is asymptotically unitary:
\begin{equation*}
\lim_{k\rightarrow\infty}\left\Vert B_{k}^{\ast}B_{k}-I\right\Vert =
\lim_{k\rightarrow\infty}\left\Vert \text{mult. by }(J_{k}^{-1}-1)\right\Vert \leq \lim_{k\rightarrow\infty}\max_{M/\!\!/G}\left\vert
J_k^{-1}-1\right\vert =0.
\end{equation*}
\
\end{proof}

\begin{theorem}
\label{thm:non-unitarity}If $\mathrm{\operatorname*{vol}}(G\cdot x_{0})$ is
not constant on $M/\!\!/G$, then there is no sequence $c_{k}$ of constants for
which $\left\Vert A_{k}A_{k}^{\ast}-c_{k}I\right\Vert $ tends to zero as $k$
tends to infinity.
\end{theorem}

In the torus case, a similar result was proved by Charles in \cite[Remark
4.29]{Charles}.

\begin{proof}
This theorem follows easily from the existence of localized sections (see
\cite[pp. 36--37]{Ma-Marinescu}). For each $[x]\in M/\!\!/G$, there exists a
sequence of holomorphic sections $\{S_{[x]}^{(k)}\in \mathcal{H}(M/\!\!/G;%
\hat{\ell}^{\otimes k})\}$ which are asymptotically concentrated near $[x]$
in the sense that for every $f\in C(M/\!\!/G)$%
\begin{equation*}
\lim_{k\rightarrow\infty}(k/2\pi)^{(n-d)/2}\int_{M/\!\!/G}f\left\vert
S_{[x]}^{(k)}\right\vert ^{2}\,\varepsilon _{\widehat{\omega}}=f([x]).
\end{equation*}
If $\mathrm{\operatorname{vol}}(G\cdot x_{0})$ is not constant on $M/\!\!/G,$
there exist points $[x_{\max }]$ and $[x_{\min }]$ where $\mathrm{\operatorname{%
vol}}(G\cdot x_{0})$ achieves its maximum and minimum values. Consider a
sequence of peak sections localized at $[x_{\max }]$. Let $\{c_{k}\}$ be a
sequence of constants. Then%
\begin{equation*}
\left\Vert A_{k}A_{k}^{\ast }-c_{k}I\right\Vert ^{2}\geq \left\Vert
A_{k}A_{k}^{\ast }S_{[x_{\max }]}^{(k)}-c_{k}S_{[x_{\max
}]}^{(k)}\right\Vert ^{2}\rightarrow \left\vert I_{k}^{-1}([x_{\max
}])-c_{k}\right\vert ^{2},~k\rightarrow \infty .
\end{equation*}%
Similarly, we obtain
\begin{equation*}
\left\Vert A_{k}A_{k}^{\ast }-c_{k}I\right\Vert ^{2}\geq \left\Vert
A_{k}A_{k}^{\ast }S_{[x_{\min }]}^{(k)}-c_{k}S_{[x_{\min
}]}^{(k)}\right\Vert ^{2}\rightarrow \left\vert I_{k}^{-1}([x_{\min
}])-c_{k}\right\vert ^{2},~k\rightarrow \infty .
\end{equation*}%
Since
\begin{equation*}
I_{k}([x_{\max }])\rightarrow \mathrm{\operatorname{vol}}(G\cdot x_{\max })\neq
\mathrm{\operatorname{vol}}(G\cdot x_{\min })\leftarrow I_{k}([x_{\min
}]),~k\rightarrow \infty ,
\end{equation*}%
we cannot have both%
\begin{equation*}
\left\vert I_{k}^{-1}([x_{\max }])-c_{k}\right\vert ^{2}\rightarrow 0~\text{%
and }\left\vert I_{k}^{-1}([x_{\min }])-c_{k}\right\vert ^{2}\rightarrow 0.
\end{equation*}%
Hence $\left\Vert A_{k}A_{k}^{\ast }-c_{k}I\right\Vert ^{2}\not\rightarrow
0. $
\end{proof}

\smallskip

The asymptotic estimates of Theorem \ref{thm:density asymptotics} can be used
to derive results about Toeplitz operators with $G$-invariant symbols. We
mention here the nicest result, concerning the asymptotics of Toeplitz
operators in the case where half-forms are included. Let $f\in C^{\infty}(M)
$, and recall that the Toeplitz operator with symbol $f$ is the map $T_{f}$
from $\mathcal{H}(M,\ell^{\otimes k}\otimes\sqrt{K})$ to itself defined by%
\[
T_{f}\,s=\operatorname{proj}(f\,s),
\]
where \textquotedblleft proj\textquotedblright\ denotes the orthogonal
projection from the space of all square-integrable sections onto the
holomorphic subspace. If $f$ is $G$-invariant, then precisely the same
argument as in the proof of Theorem \ref{thm:half-form norm decomp} shows that
for $r_{1},r_{2}\in\mathcal{H}(M,\ell^{\otimes k}\otimes\sqrt{K})$ we have%
\[
\left\langle r_{1},T_{f}r_{2}\right\rangle =\left(  k/2\pi\right)  ^{n/2}%
\int_{M}f\cdot(r_{1},r_{2})\,\varepsilon_{\omega}=\left(  k/2\pi\right)
^{(n-d)/2}\int_{M/\!\!/G}(B_{k}r_{1},B_{k}r_{2})\,J_{k}\,(f)\,\varepsilon
_{\widehat{\omega}}.
\]

If we denote by $\hat{T}_{\phi}$ the Toeplitz operator on $\mathcal{H}%
(M/\!\!/G,\hat{\ell}^{\otimes k}\otimes\sqrt{\widehat{K}})$ with symbol $\phi\in
C^{\infty}(M/\!\!/G)$, then the above formula becomes%
\[
\left\langle r_{1},T_{f}r_{2}\right\rangle =\left\langle B_{k}r_{1},\hat
{T}_{J_{k}(f)}B_{k}r_{2}\right\rangle .
\]
The asymptotics of $J_{k}(f)$ then immediately imply the following.

\begin{theorem}
\label{toeplitz.thm}Let $f\in C^{\infty}(M)^{G}$ and let $\hat{f}$ denote the
restriction of $f$ to $\Phi^{-1}(0),$ regarded as a function on $\Phi
^{-1}(0)/G.$ Let $T_{f}^{G}$ denote the restriction of the Toeplitz operator
$T_{f}$ to the space of $G$-invariant sections. Then $T_{f}^{G}$ is
asymptotically equivalent to $\hat{T}_{\hat{f}}$ in the sense that
\[
\left\Vert \hat{T}_{\hat{f}}-B_{k}\circ T_{f}^{G}\circ B_{k}^{-1}\right\Vert
\rightarrow0
\]
as $k\rightarrow\infty$.
\end{theorem}

Since we also prove that the operators $B_{k}$ are asymptotically unitary, we
may say that $T_{f}^{G}$ and $\hat{T}_{\hat{f}}$ are \textquotedblleft
asymptotically unitarily equivalent.\textquotedblright\ Such a result is
certainly does \textit{not} hold without half-forms, indicating again the
utility of the metaplectic correction.

We now turn to the proof of Theorem \ref{thm:density asymptotics}.

\begin{proof}
For all $\widehat{r},\widehat{t}\in \mathcal{H}(M/\!\!/G,\hat{\ell}^{\otimes
k}\otimes \sqrt{\widehat{K}})$, we have%
\begin{equation*}
\left\langle \widehat{r},(\widehat{T}_{\hat{f}}-B_{k}\circ T^G_{f}\circ
B_{k}^{-1})\widehat{t}\right\rangle =\left\langle \widehat{r},(\hat{f}%
-J_{k}(f))\widehat{t}\right\rangle .
\end{equation*}%
Hence, $\left\Vert \widehat{T}_{\hat{f}}-B_{k}\circ T^G_{f}\circ
B_{k}^{-1}\right\Vert =\left\Vert \text{mult. by }(\hat{f}%
-J_{k}(f))\right\Vert .$ It then follows from Theorem \ref{thm:density
asymptotics} that%
\begin{equation*}
\lim_{k\rightarrow \infty }\left\Vert \widehat{T}_{\hat{f}}-B_{k}\circ
T^G_{f}\circ B_{k}^{-1}\right\Vert =\lim_{k\rightarrow \infty }\left\Vert
\text{mult. by }(\hat{f}-J_{k}(f))\right\Vert \leq \lim_{k\rightarrow \infty
}\max_{M/\!\!/G}\left\vert \hat{f}-J_{k}(f)\right\vert =0.
\end{equation*}
\end{proof}

\smallskip

We now turn to the proof of Theorem \ref{thm:density asymptotics}. Our first task will
be to control the part of the integral near infinity. From the definition of $\tau,$
it is apparent that the integral of $\tau(\xi,x_{0})$ over $\mathfrak{g}$ is
finite for \emph{almost all} $x_{0},$ but it is not obvious that this integral
is finite for \emph{all} $x_{0}.$ We will, however, give uniform exponential
bounds on the behavior of $\tau(\xi,x_{0})$ as $\xi$ tends to infinity. This
is sufficient to show that $\tau(\xi,x_{0})\,\exp\left\{  -2k\int_{\gamma
_{\xi}}\phi_{\xi}\right\}  $ is finite for all $x_{0},$ provided $k$ is
sufficiently large, and that the part of the integral outside a neighborhood
of the origin is uniformly negligible as $k$ tends to infinity.

Our second task will be to show that Laplace's approximation can be applied so
as to give \textit{uniform} limits, which we need to prove asymptotic
unitarity of $B_{k}.$ Theorems \ref{thm:s norm decomp} and
\ref{thm:half-form norm decomp} express the densities $I_{k}$ and $J_{k}$ in
the form of (\ref{eqn:Laplace's approx}), where in the case of the expression
for $J_{k},$ the $k$-independent term in the exponent should be factored out
and grouped with $\tau(\xi,x_{0}).$ We develop a uniform version of Laplace's
approximation that allows for the desired uniformity. (It is probably well
known that this is possible, but we were not able to find a written proof.)
Our argument will be based on the Morse--Bott Lemma (a parameterized version
of the usual Morse lemma).

Once these two tasks are accomplished, it is a straightforward matter to plug
in the Hessian computation in Theorem \ref{thm:s along a curve} to obtain
Theorem \ref{thm:density asymptotics}.

\subsection{Growth estimates}

In this section, we show that the contribution to $I_{k}$ coming from the
complement of a tubular neighborhood of the zero-set is negligible as $k$
tends to infinity:

\begin{theorem}
\label{lemma:negligable contributions}There exist constants $\,b,D>0$ such
that for all $x_{0}\in\Phi^{-1}(0)$ and for all $R$ and $k$ sufficiently
large
\[
\int_{\mathfrak{g}\setminus B_{R}(0)}f(\xi,x_{0})\tau(\xi,x_{0})\,\exp\left\{
-\int_{\gamma_{\xi}}2k\,\phi_{\xi}\right\}  \,d^{d}\xi\leq be^{-RDk},
\]
where $B_R(0)$ is a ball of radius $R$ centered at $0\in\mathfrak{g}.$

\end{theorem}

Our first step is to show that the integrand appearing in $I_{k}$ decays
exponentially in the radial directions.

\begin{lemma}
\label{lemma:exp decay}There exists $C>0$ such that for all sufficiently large
$t$,
\[
\exp\left\{  -\int_{\gamma_{t\hat{\xi}}}2k\,\phi_{t\hat{\xi}}\right\}  \leq
e^{-2ktC}%
\]
uniformly on $\Phi^{-1}(0),$ where $\hat{\xi}\in\mathfrak{g}$ with $\left\vert
\hat{\xi}\right\vert =1$.
\end{lemma}

\begin{proof}
By definition, we have%
\begin{equation*}
\int_{\gamma_{t\hat{\xi}}}\phi_{t\hat{\xi}}=\int_{0}^{1}\langle\Phi(e^{i\tau
t\hat{\xi}}\cdot x_{0}),t\hat{\xi}\rangle\,d\tau=t\int_{0}^{1}\langle
\Phi(e^{i\tau t\hat{\xi}}\cdot x_{0}),\hat{\xi}\rangle\,d\tau.
\end{equation*}
Hence, it is enough to show that there exists some $C>0$ such that for $t$
sufficiently large and for all $x_{0}\in\Phi^{-1}(0)$%
\begin{equation*}
\int_{0}^{1}\langle\Phi(e^{i\tau t\hat{\xi}}\cdot x_{0}),\hat{\xi}%
\rangle\,d\tau\geq C.
\end{equation*}
For $t\geq 0,$ define
\begin{equation}
M_{t}=\left\{ e^{it\xi }\cdot x_{0}:\left\vert \hat{\xi}\right\vert
=1,~x_{0}\in \Phi ^{-1}(0)\right\} .  \label{eqn:M_t}
\end{equation}%
Then $M_{t},~t>0,$ is a sphere bundle over $\Phi ^{-1}(0)$ and is hence
compact. Define $\rho _{t}:M_{t}\rightarrow \mathbb{R}$ by
\begin{equation*}
\rho _{t}(e^{it\hat{\xi} }\cdot x_{0})=2\int_{0}^{1}\phi _{\hat{\xi}}(e^{i\tau t%
\hat{\xi}}\cdot x_{0})\,d\tau .
\end{equation*}%
Fix $t>0.$ Then since $\rho _{t}$ is smooth with compact domain there is
some $m\in M_{t}$ where $\rho _{t}$ achieves its minimum.
Recall that $e^{it\tau \hat{\xi}}\cdot x_{0}$ is the gradient line of $\phi
_{t\hat{\xi}}$. Moreover,%
\begin{equation*}
(\operatorname{grad}\phi _{\hat{\xi}})_{e^{i\xi }\cdot x_{0}}=JX_{e^{i\xi }x_{0}}^{%
\hat{\xi}}\text{ and }\left\vert JX_{e^{i\xi }\cdot x_{0}}^{\hat{\xi}%
}\right\vert =\left\vert X_{e^{i\xi }\cdot x_{0}}^{\hat{\xi}}\right\vert
\neq 0
\end{equation*}%
since $G$ acts freely on $M_{s}.$ Therefore, $\phi _{t\hat{\xi}}$ is
strictly increasing along $e^{it\tau \hat{\xi}}\cdot x_{0}$. In particular, $%
2C:=\rho _{t}(m)>0$ so that, by the definition of $m,$ we have $\rho
_{t}(m^{\prime })\geq 2C$ for all $m^{\prime }\in M_{t}.$ Hence we obtain $%
\int_{0}^{1}\langle \Phi (e^{i\tau t^{\prime }\hat{\xi}}\cdot x_{0}),t\hat{%
\xi}\rangle \,d\tau \geq C$ for all $t^{\prime }>t$ and for all $x_{0}\in
\Phi ^{-1}(0).$
\end{proof}

We now turn to the task of estimating $\tau,$ which we will do by embedding
$M$ into projective space. We recall some basic projective geometry. The
function $\tau$ is the Jacobian of the map $\psi:\mathfrak{g}\times\Phi
^{-1}(0)\rightarrow M.$ Since $M$ is a compact K\"{a}hler manifold and $\ell$
is an equivariant ample line bundle, we can equivariantly embed $M$ into
$\mathbb{C}P^{N}$ for some $N$ (this is Kodaira's embedding theorem enhanced
by the presence of a $G$-invariant K\"{a}hler structure and line bundle; the
embedding is realized by considering holomorphic sections of
high tensor powers of $\ell$). We identify $M$ with its embedded image. Under
this embedding, we identify $G$ with a compact subgroup of $PU(N+1)$.

The metric on the embedded image of $M$ induced by the K\"{a}hler metric on
$M$ differs from the Fubini--Study metric on $\mathbb{C}P^{N}$ by smooth
factors. In particular, the difference between computing a determinant with
the Fubini--Study metric and the induced K\"{a}hler metric is the determinant
of the linear map which takes a Fubini--Study orthogonal basis of the tangent
space to a basis which is orthogonal with respect to the induced metric. This
map is invertible and smooth on $M$ and hence is bounded and bounded away from
zero, and so it suffices for our purposes to compute the Jacobian $\tau$ using
the Fubini--Study metric.

We use homogeneous coordinates on $\mathbb{C}P^{N};$ i.e.,$\mathbb{C}%
P^{N}=\{[\mathbf{z}=(z^{0},z^{1},\dots,z^{N})]:z^{j}\in\mathbb{C}\}.$ The
tangent space $T_{[\mathbf{z}]}\mathbb{C}P^{N}$ can be identified with the
orthogonal complement of the line in $\mathbb{C}^{N+1}$ determined by
$\mathbf{z}$.

The Fubini--Study metric on $\mathbb{C}P^{N}$ is \cite[Sec. 7.4]{Zheng}%
\[
g_{\text{FS}}=\operatorname{Re}\left(  \frac{\left\vert \mathbf{z}\right\vert
^{2}d\mathbf{\bar{z}}\otimes d\mathbf{z-}\sum_{j,k=0}^{N}z^{j}d\bar{z}%
^{j}\otimes\bar{z}^{k}dz^{k}}{\left\vert \mathbf{z}\right\vert ^{4}}\right)  .
\]

Now we show that $\tau(\xi,x_{0})$ grows at most exponentially in the radial
directions away from the zero-set.

\begin{lemma}
\label{lemma:tau estimates}There exist constants $a,b>0$ such that for all $t
$ sufficiently large
\[
\tau(t\hat{\xi},x_{0})\leq bt^{-d}e^{at}%
\]
uniformly on $\Phi^{-1}(0)$ and for all $\hat{\xi}\in\mathfrak{g}$ with
$\left\vert \hat{\xi}\right\vert =1.$
\end{lemma}

\begin{proof}
Recall that $\Xi=$ $\{\xi_{j}\}_{j=1}^{d}$ is an orthonormal basis of $%
\mathfrak{g}<\mathfrak{pu}(N+1)$ with respect to an $\operatorname{Ad}$-invariant
inner product on $\mathfrak{g}$. Let $\{X_{j}\}_{j=1}^{2n-d}$ and $%
\{Y_{j}\}_{j=1}^{2n}$be orthonormal bases of $T_{x_{0}}\Phi^{-1}(0)$ and $%
T_{e^{it\hat{\xi}}\cdot x_{0}}M$, respectively. Write a vector in terms of
its components with respect to the appropriate basis using the usual
convention; for example, $v=v^{j}Y_{j}\in T_{e^{i\xi}\cdot x_{0}}M$. We want
to estimate%
\begin{equation*}
\tau(t\hat{\xi},x_{0})=(\det\psi_{\ast})(t\hat{\xi},x_{0})
\end{equation*}
We can compute the components of the images of the basis $\Xi \cup\{X_{j}\}$
by
\begin{equation}
\psi_{\ast}(\xi_{j})^{l}=g_{\text{FS}}(Y_{l},\psi_{\ast}(\xi_{j})),~\text{%
and }\psi_{\ast}(X_{j})^{l}=g_{\text{FS}}(Y_{l},\psi_{\ast}(X_{j})).
\label{eqn:components}
\end{equation}
There is a unique $R\in\mathfrak{su}(N+1)$ such that $e^{i\hat{\xi}} [%
\mathbf{z}]=[e^{iR}\mathbf{z}].$ Since $R$ is skew-Hermitian, it has pure
imaginary eigenvalues and these eigenvalues vary continuously with respect
to the choice of $\xi$. Let the eigenvalues be $-\sqrt{-1}\lambda_{0},\dots,-%
\sqrt{-1}\lambda_{N}$ with $\lambda_{0}\leq\cdots\leq\lambda_{N}.$ Then for
generic $\mathbf{z}\in\mathbb{C}^{N+1},$ we have $\left\vert (e^{itR}\mathbf{%
z)}^{j}\right\vert =O(e^{\lambda_{N}t}).$
In order to compute the components in equation (\ref{eqn:components}), we
need to estimate the growth of both $Y_{l}$ and the pushforwards $%
\psi_{\ast}(\xi_{j})$ and $\psi_{\ast}(X_{j})$. We begin with $Y_{l}$.
Let $[\mathbf{w}_{t}]=e^{it\hat{\xi}}\cdot x_{0}.$ Since $Y_{j}$ $%
=(Y_{j}^{0},\dots,Y_{j}^{N})$ is a unit vector, we have%
\begin{equation*}
\left( \left\vert Y_{j}\right\vert _{\text{FS}}^{2}\right)
_{e^{it\xi}x_{0}}=1=\frac{1}{\left\vert \mathbf{w}_{t}\right\vert ^{4}}\operatorname{%
Re}\left( \left\vert \mathbf{w}_{t}\right\vert ^{2}Y_{j}\cdot\bar{Y}%
_{j}-\sum_{l,m=0}^{N}w_{t}^{l}\bar{w}_{t}^{m}\,\bar{Y}_{j}^{l}Y_{j}^{m}%
\right) .
\end{equation*}
In particular, since $\left\vert \mathbf{w}_{t}\right\vert
^{4}=O(e^{4\lambda _{N}t})$ and $\left\vert w_{t}^{l}\right\vert
=O(e^{\lambda_{N}t})$, we must have $\left\vert Y_{j}^{l}\right\vert
=O(e^{\lambda_{N}t}).$
To estimate the growth rate of the pushforward $\psi_{\ast}(X_{j})$, we
observe that if $X_{j}=[\mathbf{x}_{j}]$, then
\begin{equation*}
\psi_{\ast}(X_{j})=[e^{itR}\mathbf{x}_{j}]
\end{equation*}
from which we obtain%
\begin{equation}
\psi_{\ast}(X_{j})^{l}=\frac{1}{\left\vert \mathbf{w}_{t}\right\vert ^{4}}%
\operatorname{Re}\left( \left\vert \mathbf{w}_{t}\right\vert ^{2}\sum _{m}\bar{Y}%
_{l}^{m}(e^{itR}\mathbf{x}_{j})^{m}-\sum_{l,m}w_{t}^{l}\bar{w}_{t}^{m}\,\bar{%
Y}_{l}^{l}(e^{itR}\mathbf{x}_{j})^{m}\right) .  \label{eqn:pushforward1}
\end{equation}
By the same argument as above, we see that $\left\vert (e^{itR}\mathbf{x}%
_{j})^{m}\right\vert =O(e^{\lambda_{N}t}).$
We must be a little careful with our growth estimates; if $\lambda_{0}>0,$
then at worst
\begin{equation*}
\left\vert \mathbf{w}_{t}\right\vert ^{-1}=O(e^{-\lambda_{0}t})=O(1),
\end{equation*}
but if $\lambda_{0}<0$, then there may be certain points $x_{0}$ where $%
\left\vert \mathbf{w}_{t}\right\vert ^{-1}=O(e^{-\lambda_{0}t}),$ which is
exponential \emph{growth}. With this in mind, we can use equation (\ref%
{eqn:pushforward1}) to make our first growth estimates:%
\begin{equation*}
\left\vert \psi_{\ast}(X_{j})_{e^{it\hat{\xi}}\cdot x_{0}}^{l}\right\vert
=O(e^{4(\lambda_{N}-\lambda_{0})t}).
\end{equation*}
Of course, it is clear from equation (\ref{eqn:pushforward1}) that at
generic points in generic directions, the growth will actually be much less
than this estimate indicates since the numerator and denominator will both
grow at the same rate.
Finally we consider the components $\psi_{\ast}(\xi_{j})^{l}$. Since we have
embedded $M$ in $\mathbb{C}P^{N},$ the exponential $e^{itR}$ is a matrix
exponential, and there is an explicit formula for its derivative. Write $%
x_{0}=[\mathbf{z}]$ and let $S_{j}\in\mathfrak{su}(N+1)$ be the matrix such
that $e^{i(t\hat{\xi}+s\xi_{j})}\cdot x_{0}=[e^{i(tR+sS_{j})}\mathbf{z}]$.
Using the formula for the derivative of the exponential, we compute%
\begin{equation*}
\psi_{\ast}(\xi_{j})=\frac{d}{ds}|_{s=0}e^{i(t\hat{\xi}+s\xi_{j})}\cdot
x_{0}=\left[ e^{itR}\left( \frac{1-e^{-it\operatorname{ad}(R)}}{t\operatorname{ad}(R)}%
\right) S_{j}\mathbf{z}\right] .
\end{equation*}
As a linear operator on $\mathfrak{su}(N+1)$, the matrix $\operatorname{ad}(R)$
is skew and hence has purely imaginary eigenvalues (which vary continuously
with respect to choice of $\xi$). Denote them by $\sqrt{-1}\mu_{0},\dots,%
\sqrt{-1}\mu_{d}$ with $\mu_{0}\leq\cdots\leq\mu_{d}.$ Then%
\begin{equation*}
\left\vert \left( e^{itR}\left( \frac{1-e^{-it\operatorname{ad}(R)}}{t\operatorname{ad}%
(R)}\right) S_{j}\mathbf{z}\right) ^{m}\right\vert =O\left( \frac{%
e^{(\lambda_{N}+\mu_{d})t}}{t}\right) .
\end{equation*}
Finally, we can use this to estimate the growth of the component $\psi_{\ast
}(\xi_{j})^{l}$ by computing with the Fubini--Study metric as above:%
\begin{align*}
\psi_{\ast}(\xi_{j})^{l} & =\frac{1}{\left\vert \mathbf{w}_{t}\right\vert
^{4}}\Big(\left\vert \mathbf{w}_{t}\right\vert ^{2}\sum_{m}\bar{Y}%
_{l}^{m}\left( e^{itR}\left( \frac{1-e^{-it\operatorname{ad}(R)}}{t\operatorname{ad}(R)%
}\right) S_{j}\mathbf{z}\right) ^{m} \\
& \qquad-\sum_{l,m}w_{t}^{l}\bar{Y}_{l}^{m}\,\bar{w}_{t}^{m}\left(
e^{itR}\left( \frac{1-e^{-it\operatorname{ad}(R)}}{t\operatorname{ad}(R)}\right) S_{j}%
\mathbf{z}\right) ^{m}\Big) \\
& =O\left( \frac{e^{\left( 4\lambda_{0}+3\lambda_{N}+\mu_{d}\right) t}}{t}%
\right) .
\end{align*}
Putting this all together, we obtain an estimate of the growth of the
Jacobian $\tau$. Define%
\begin{equation*}
a(\hat{\xi})=d(4\lambda_{0}+3\lambda_{N}+\mu_{d})+4(2n-d)(\lambda_{N}-%
\lambda_{0}).
\end{equation*}
Then for all $x_{0}\in\Phi^{-1}(0)$%
\begin{align*}
\left\vert \tau(t\hat{\xi},x_{0})\right\vert & =O\left( \left(
t^{-1}e^{(4\lambda_{0}+3\lambda_{N}+\mu_{d})t}\right) ^{d}\left(
e^{4(\lambda_{N}-\lambda_{0})t}\right) ^{2n-d}\right) \\
& =O\left( t^{-d}e^{(d(4\lambda_{0}+3\lambda_{N}+\mu_{d})+4(2n-d)(\lambda
_{N}-\lambda_{0}))t}\right) =O\left( t^{-d}e^{a(\hat{\xi})t}\right) .
\end{align*}
Since the eigenvalues $\lambda_{j}$ and $\mu_{j}$ depend continuously on the
choice of $\xi$, and the unit sphere $\{\hat{\xi}\in\mathfrak{g}:\hat{\xi }%
=1\}$ is compact, there is some point $\hat{\xi}_{\max}\in\mathfrak{g}$
where $a(\hat{\xi})$ attains its maximum; let $a=a(\hat{\xi}_{\max}).$ Then $%
\left\vert \tau(t\hat{\xi},x_{0})\right\vert =O(t^{-d}e^{at})$ uniformly in $%
x_{0}$ and $\hat{\xi}$ as desired.
\end{proof}

Next, we combine the previous two lemmas and integrate over the complement of
a neighborhood of $0\in\mathfrak{g}$ to prove Theorem
\ref{lemma:negligable contributions}.

\begin{proof}[Proof of Theorem \protect\ref{lemma:negligable contributions}]
Introduce polar coordinate $\xi =(t,\hat{\xi})$, $t\in \mathbb{R}_{+}$ on $%
\mathfrak{g}$. Then $d^{d}\xi =t^{d-1}J(\hat{\xi})dt\wedge d\hat{\xi}$,
where $d\hat{\xi}$ is the solid angle element and $J(\hat{\xi})$ is the
appropriate Jacobian. By the previous two lemmas and the fact that $f\in
C^{\infty }(M)$ implies $f$ is bounded on $M$, there is some $b^{\prime }>0$
such that for $R$ and $k$ sufficiently large%
\begin{align*}
\Bigg\vert\int_{\mathfrak{g}\setminus B_{R}(0)}f\,\tau \,& \exp \left\{
-\int_{\gamma _{\xi }}2k\,\phi _{\xi }\right\} \,d^{d}\xi \Bigg\vert \\
& =\left\vert \int_{R}^{\infty }\int_{S^{d-1}}\tau (t\hat{\xi},x_{0})\exp
\left\{ -\int_{\gamma _{t\hat{\xi}}}2k\,\phi _{t\hat{\xi}}\right\} t^{d-1}J(%
\hat{\xi})d\hat{\xi}dt\right\vert  \\
& \leq b^{\prime }\int_{R}^{\infty
}\int_{S^{d-1}}t^{-d}e^{at}e^{-2kCt}t^{d-1}J(\hat{\xi})\,d\hat{\xi}dt.
\end{align*}%
For $k$ sufficiently large there is some $b^{\prime \prime }>0$ such that%
\begin{equation*}
\int_{R}^{\infty }\int_{S^{d-1}}\frac{e^{(a-2kC)t}}{t}J(\hat{\xi})\,dt\wedge
d\hat{\xi}=\mathrm{\operatorname{vol}}(S^{d-1})\int_{R}^{\infty }\frac{%
e^{(a-2kC)t}}{t}dt\leq b^{\prime \prime }e^{-RDk}.
\end{equation*}%
Letting $b=b^{\prime }b^{\prime \prime }$, we obtain our desired result.
\end{proof}

\subsection{Laplace's approximation}

We now turn to the issue of uniformity in Laplace's approximation. We need to
show that when applied to the expressions for $I_{k}(f)([x_{0}])$ and
$J_{k}(f)([x_{0}]),$ we obtain limits that are uniform in $x_{0}.$ To do this,
follow a standard argument for Laplace's approximation using the Morse Lemma,
but we use a parameterized version (the Morse--Bott Lemma). Then if we simply
are careful to bound our errors at each stage, we will obtain the desired
uniform limits.

We now state the Morse--Bott Lemma.\ (A careful proof can be found in
\cite{Banyaga-Hurtubise}.) A point $p\in M$ is a critical point of a function
$\rho:M\rightarrow\mathbb{R}$ if the differential of $\rho$ at $p$ is zero. A
function $\rho$ is said to be a Morse--Bott function if $Crit(\rho) $, the set
of critical points of $\rho$, is a disjoint union of connected smooth
submanifolds and for each critical point, the Hessian $H(\rho)(p)$, in the
directions normal to $C$, is non-degenerate.

\begin{lemma}
[Morse--Bott Lemma]Let $\rho:M\rightarrow\mathbb{R}$ be a Morse--Bott
function, $C$ a connected component of $Crit(f)$ of dimension $2n-d$, and
$p\in C$. Then there exists and open neighborhood $U$ of $p$ and a smooth
chart $\varphi:U\rightarrow\mathbb{R}^{d}\times\mathbb{R}^{2n-d}$ such that

\begin{enumerate}
\item $\varphi(p)=0$

\item $\varphi(U\cap C)=\{(x,y)\in\mathbb{R}^{d}\times\mathbb{R}^{2n-d}:x=0\}
$ and

\item $(\rho\circ\varphi^{-1})(z,y)=\rho(C)-\frac{1}{2}(z_{1}^{2}+z_{2}%
^{2}+\cdots+z_{k}^{2})+\frac{1}{2}(z_{k+1}^{2}+\cdots+z_{2n-d}^{2})$

where $k\leq2n-d$ is the index of $H(\rho)(p)$ and $f(C)$ is the common value
of $\rho$ on $C$. Moreover, if $\varphi^{\prime}$ is another smooth chart
satisfying (1) and (2) near $p$, the Jacobian of the coordinate change in the
$\mathbf{z}$-directions, evaluated at $p$, is $\det\frac{dz}{dz^{\prime}%
}=\sqrt{\det H(\rho)(p)}$.
\end{enumerate}
\end{lemma}

To apply Laplace's approximation in our case, we need to know the value of the
density $\tau(\xi,x_{0})$ when $\xi=0.$

\begin{lemma}
\label{lemma:tau=1}The function $\tau$ equals $\mathrm{\operatorname*{vol}%
}(G\cdot x_{0})$ on the zero-set of the moment map.
\end{lemma}

\begin{proof}
Let $W_{x}=\{X_{x}^{\xi}:\xi\in\mathfrak{g}\}$. Then we claim that the
complexified tangent space to $M,$ at a point $x_{0}\in\Phi^{-1}(0),$
decomposes as a $B$-orthogonal direct sum
\begin{equation}
T_{x_{0}}^{\mathbb{C}}M=W_{x_{0}}^{\mathbb{C}}\oplus JW_{x_{0}}^{\mathbb{C}%
}\oplus(T_{x_{0}}^{1,0}M\cap T_{x_{0}}\Phi^{-1}(0)).  \label{eqn:orthogonal}
\end{equation}
To see this, first note that by the definition of the $G_{\mathbb{C}}$%
-action, the tangent space to the $G_{\mathbb{C}}$-orbit $G_{\mathbb{C}%
}\cdot x$ is $W_{x}\oplus JW_{x}.$ Since $G\cdot x_{0}$ is a totally real
submanifold,
\begin{equation*}
W_{x_{0}}\cap T_{x_{0}}^{1,0}M=\{0\},
\end{equation*}
and so a dimension count shows that the complexified tangent space
decomposes into the claimed direct sum.

We need to show that the direct sums are in fact $B$-orthogonal. First, if $%
v\in T_{x_{0}}\Phi^{-1}(0)$, then $v(\phi_{\xi})=0$ since the moment map is
constant (and in fact identically zero) on $\Phi^{-1}(0).$ Hence $%
B(v,JX^{\xi })=\omega(v,X^{\xi})=-d\phi_{\xi}(v)=-v(\phi_{\xi})=0$ and so $%
JW_{x_{0}}$ is $B$-orthogonal to $T_{x_{0}}\Phi^{-1}(0).$
Next, suppose $v\in(T_{x_{0}}^{1,0}M\cap T_{x_{0}}\Phi^{-1}(0)).$ Then $%
Jv=i\,v$ so that
\begin{equation*}
B(X^{\xi}+JX^{\xi^{\prime}},v)=i\,\omega(X^{\xi},v)+i\omega(JX^{\xi^{%
\prime}},v)=iv(\phi_{\xi})+v(\phi_{\xi})=0
\end{equation*}
since $v$ is a linear combination of
vectors in $T_{x_{0}}\Phi^{-1}(0)$. This shows $W_{x_{0}}^{\mathbb{C}%
}\oplus JW_{x_{0}}^{\mathbb{C}}$ is $B$-orthogonal to $(T_{x_{0}}^{1,0}M\cap
T_{x_{0}}\Phi^{-1}(0))$.
Finally, $B(X^{\xi},JX^{\xi^{\prime}})=\omega(X^{\xi},X^{\xi^{\prime}})=\{%
\phi_{\xi},\phi_{\xi^{\prime}}\}=\phi_{\lbrack\xi,\xi^{\prime}]}=0$ on the
zero-set, so $W_{x_{0}}$ is $B$-orthogonal to $JX_{x_{0}}$.

Having established (\ref{eqn:orthogonal}), we now turn to the computation of
$\tau.$ On the zero-set of the moment map, the tangent space to the orbit $%
\exp(i\mathfrak{g})\cdot x_{0}$ is $JW_{x_{0}},$ which, by (\ref%
{eqn:orthogonal}), is $B$-orthogonal to the tangent space of the zero-set.
So if we choose coordinates $\{x^{j}\}$ on $M$ near $x_{0}\in \Phi^{-1}(0)$
such that $\exp(i\mathfrak{g})\cdot x_{0}=\{x^{d+1}=\cdots=x^{2n}=0\}$, and $%
x^{j}(e^{i\xi}\cdot x_{0})=\langle\xi_{j},\xi\rangle$ (i.e.,$x^{j}$ is the
image in $M$ of the linear coordinate in the direction of $\xi_{j}$ on $%
\mathfrak{g}$), then $B$ is block diagonal and so for $x\in \Phi^{-1}(0)$
\begin{equation*}
\det B_{x}=(\det B_{x})|_{\exp(i\mathfrak{g})\cdot x_{0}}(\det
B_{x})|_{\Phi^{-1}(0)}.
\end{equation*}
By the standard formula for the volume form in local coordinates, we have%
\begin{align*}
(\varepsilon_{\omega})_{x} & =d\mathrm{\operatorname{vol}}(M)_{x}=\sqrt{\det B_{x}%
}dx^{1}\wedge\cdots\wedge dx^{2n} \\
& =\sqrt{\det B_{x}|_{\exp(i\mathfrak{g})\cdot x_{0}}}dx^{1}
\wedge\cdots\wedge dx^{d} \wedge \sqrt{\det B_{x}|_{\Phi ^{-1}(0)}}%
dx^{d+1}\wedge\cdots\wedge dx^{2n} \\
& =\sqrt{\det B_{x}|_{\exp(i\mathfrak{g})\cdot x_{0}}}d^{d}\xi\wedge \,d%
\mathrm{\operatorname{vol}}(\Phi^{-1}(0))_{x_{0}}\,.
\end{align*}
We can rewrite $\det B_{x}|_{\exp(i\mathfrak{g})\cdot x_{0}}=\det\left(
B\left( \frac{\partial}{\partial x^{j}},\frac{\partial}{\partial x^{k}}%
\right) \right) $ in the desired form by noting that $\frac{\partial }{%
\partial x^{j}}=\frac{d}{dt}|_{t=0}e^{it\xi_{j}}x=JX_{x}^{\xi_{j}}$ so that
\begin{equation*}
\det\left( B\left( \frac{\partial}{\partial x^{j}},\frac{\partial}{\partial
x^{k}}\right) \right) =\det\left( B\left( JX^{\xi_{j}},JX^{\xi_{k}}\right)
\right) =\det\left( B\left( X^{\xi_{j}},X^{\xi_{k}}\right) \right)
=\det\nolimits_{\Xi}B.
\end{equation*}
The lemma now follows from Lemma \ref{lemma:det = vol}.
\end{proof}

We are now ready to perform a uniform Laplace's approximation. We follow a
standard argument (see, for example, \cite[Sec. 8.3]{Bleistein}), although we
keep careful track of the error terms to prove uniformity of the limit.

Denote the ball of radius $R>0$ centered at $0\in\mathfrak{g}$ by $B_{R}(0)$.

\begin{lemma}
\label{lemma:density limit}Let $f\in C^{\infty}(M)$ and define%
\[
I_{k,R}(f)(x_{0})=\left(  k/2\pi\right)  ^{d/2}\mathrm{\operatorname*{vol}%
}(G\cdot x_{0})\int_{B_{R}(0)}f(\xi,x_{0})e^{-k\rho(\xi,x_{0})}\tau(\xi
,x_{0})\,d^{d}\xi
\]
where $\rho(\xi,x_{0})=2\int_{0}^{1}\phi_{\xi}(e^{it\xi}\cdot x_{0})\,dt.$
Then there exists some $R>0$ such that
\[
\lim_{k\rightarrow\infty}\left\vert I_{k,R}(f)(x_{0})-2^{-d/2}%
\mathrm{\operatorname*{vol}}(G\cdot x_{0})f(x_{0})\right\vert =0
\]
uniformly on $\Phi^{-1}(0),$ where $H(\rho)$ is the Hessian of $\rho.$
\end{lemma}

\begin{proof}
Laplace's method is based on an application of the Morse--Bott lemma; the
zero-set $\{0\}\times \Phi ^{-1}(0)$ is a critical submanifold for the
function $\rho $---in fact, $\rho (0,\Phi ^{-1}(0))=0$ is a minimum since $%
e^{it\xi }x_{0}$ is the gradient flow line of $\phi _{\xi }$---and so there
exist local coordinates $(\mathbf{z},\mathbf{y})$ in some neighborhood of $x_{0}$ such that
\begin{equation}
(\rho \circ \varphi ^{-1})(\mathbf{z},\mathbf{y})=\frac{1}{2}\mathbf{z}\cdot
\mathbf{z}.  \label{eqn:morse coords}
\end{equation}%
\newline

For each $x_{0}$ the coordinates $(\mathbf{z,y})$ are defined in some ball
of positive radius; since the zero-set $\{0\}\times \Phi ^{-1}(0)$ is
compact we can define $R>0$ to be the minimum of these radii.
If we regard the coordinates $\mathbf{z}$ as new coordinates on $\mathfrak{g}
$ then there is a Jacobian $J(\mathbf{z})=\det \frac{\partial \mathbf{\xi }}{%
\partial \mathbf{z}}$ of the coordinate change in the $\mathfrak{g}$%
-directions, whence%
\begin{equation*}
I_{k,R}(f)(x_{0})=\left( k/2\pi \right) ^{d/2}\mathrm{\operatorname{vol}}(G\cdot
x_{0})\int_{B_{R}(0)}e^{-\frac{k}{2}\mathbf{z}^{2}}f(\mathbf{z},x_{0})\tau (%
\mathbf{z},x_{0})J(\mathbf{z})\,d^{d}\mathbf{z}.
\end{equation*}

The exact first-order Taylor series of the function $T(\mathbf{z},x_{0})=f(%
\mathbf{z},x_{0})\tau (\mathbf{z},x_{0})J(\mathbf{z})$ is given by expanding
the identity $T(\mathbf{z},x_{0})=\int_{0}^{1}\frac{d}{dt}T(t\mathbf{z}%
,x_{0})dt$ to yield $T(\mathbf{z},x_{0})=T(0,x_{0})+\mathbf{z}\cdot \mathbf{S%
}(\mathbf{z},x_{0}),$ where $S_{j}(\mathbf{z},x_{0})=\int_{0}^{1}\frac{%
\partial }{\partial z^{j}}T(t\mathbf{z},x_{0})dt.$
We now have%
\begin{equation*}
(2\pi/k)^{d/2}\frac{I_{k,R}(x_{0})} {\mathrm{\operatorname{vol}}(G\cdot x_{0})}%
=T(0,x_{0})\int_{B_{R}(0)}e^{-\frac{k}{2}\mathbf{z}^{2}}d^{d}\mathbf{z}%
+\int_{B_{R}(0)}e^{-\frac {k}{2}\mathbf{z}^{2}}(\mathbf{z}\cdot\mathbf{S}%
)d^{d}\mathbf{z.}
\end{equation*}

Let $\vec{\nabla}=\left( \frac{\partial}{\partial z^{1}},\dots,\frac {%
\partial}{\partial z^{d}}\right)$ be the vector calculus gradient. Then
\begin{equation*}
-\frac{k}{2}e^{-\frac{k}{2}\mathbf{z}^{2}}(\mathbf{z}\cdot\mathbf{S} )=\vec{%
\nabla}\cdot(\mathbf{S}e^{-\frac{k}{2}\mathbf{z}^{2}})-(\vec{\nabla }\cdot%
\mathbf{S})e^{-\frac{k}{2}\mathbf{z}^{2}}.
\end{equation*}
Substituting this into the last term above and applying the divergence
theorem yields
\begin{equation*}
\int_{B_{R}(0)}\vec{\nabla}\cdot(\mathbf{S}e^{-\frac{k}{2}\mathbf{z}%
^{2}})d^{d}\mathbf{z} = \int_{\partial B_{R}(0)}(\mathbf{\hat{n}\cdot
S})e^{-\frac{k}{2}\mathbf{z}^{2}}d^{d-1}\mathbf{\Omega}\leq e^{-\frac{k}{2}%
r^{2}}\int_{\partial B_{R}(0)}(\mathbf{\hat{n}\cdot S})\,d^{d}\mathbf{z}
\end{equation*}
where $\mathbf{\hat{n}}$ is an outward pointing unit normal vector to $%
\partial B_{R}(0)$ and $d^{d-1}\mathbf{\Omega}$ is the solid angle element. The
last integral above is a continuous function of $x_{0}\in\Phi^{-1}(0)$, so it is
bounded by a constant $Q_{1}>0$. Similarly, there is some $Q_{2}>0$ such that
\[Q_{2}=\max_{\underset{x_{0}\in\Phi^{-1}(0)}{\mathbf{z}\in B_R(0)}}\left\vert\vec{\nabla}\cdot\mathbf{S}\right\vert .\]
We can now estimate that for any $x_{0}\in\Phi^{-1}(0)$%
\begin{equation*}
\Bigg\vert (2\pi/k)^{d/2}\frac{I_{k,R}(x_{0})}{\mathrm{\operatorname{vol}}(G\cdot
x_{0})}-T(0,x_{0}) \int_{B_{R}(0)}e^{-\frac{k}{2}\mathbf{z}^{2}}d^{d}\mathbf{%
z}\Bigg\vert\leq\frac{Q_{2}}{k}\int_{B_{R}(0)}e^{-\frac{k}{2}\mathbf{z}%
^{2}}d^{d}\mathbf{z+}\frac{Q_{1}}{k}e^{-\frac{k}{2}r^{2}}.
\end{equation*}

It is a standard result \cite[Lemma 8.3.1]{Bleistein} that for any bounded
domain $D$ and for every $m>0$, there exists a constant $C_{m}>0$ such that%
\begin{equation*}
\int_{D}e^{-\frac{k}{2}\mathbf{z}^{2}}d^{d}\mathbf{z\leq\,}\left( \frac{2\pi
}{k}\right) ^{d/2}+C_{m}\lambda^{-m}.
\end{equation*}
Putting this all together, we have for every $x_{0}\in\Phi^{-1}(0)$ and for
every $m>0$%
\begin{align*}
\Bigg\vert (2\pi/k)^{d/2}\frac{I_{k,R}(x_{0})}{\mathrm{\operatorname{vol}}(G\cdot
x_{0})} & -T(0,x_{0})\Bigg\vert \\
& \leq\left\vert (2\pi/k)^{d/2}I_{k,R}(x_{0})-T(0,x_{0})\int_{B_{R}(0)}e^{-%
\frac{k}{2}\mathbf{z}^{2}}d^{d}\mathbf{z}\right\vert \\
& \qquad\qquad+\left\vert T(0,x_{0})\right\vert \left\vert
\int_{B_{R}(0)}e^{-\frac{k}{2}\mathbf{z}^{2}}d^{d}\mathbf{z}%
-(2\pi/k)^{d/2}\right\vert \\
& \leq\frac{Q_{2}}{k}(2\pi/k)^{d/2}+\frac{Q_{1}}{k}e^{-\frac{k}{2}%
r^{2}}+(Q_{2}+Q_{3})C_{m}k^{-m}
\end{align*}
where $Q_{3}=\max_{x_{0}\in\Phi^{-1}(0)}\left\vert T(0,x_{0})\right\vert .$

All that remains is to note that by the Morse--Bott lemma and Lemma \ref%
{lemma:tau=1},
\[T(0,x_{0})=\frac{f(x_{0})\mathrm{\operatorname{vol}}(G\cdot x_{0})%
}{\sqrt{\det H(\rho )(0,x_{0})}}.\]
By Theorem \ref{thm:s along a curve} and Lemma \ref{lemma:det = vol}, the
denominator is $2^{d/2}\mathrm{\operatorname{vol}}(G\cdot x_{0})$. Putting this
into the above inequality and taking the limit $k\rightarrow \infty$ we obtain
our desired result.
\end{proof}

If we include the metaplectic correction, then a similar proof applies if we
replace $\tau(\xi,x_{0})$ by $\tau(\xi,x_{0})\,\exp\left\{  -\int_{\gamma
_{\xi}}\frac{{\mathcal{L}}_{JX^{\xi}}\varepsilon_{\omega}}{2\varepsilon
_{\omega}}\right\}  .$ Since the argument of the exponent is zero on the
zero-set, we obtain a similar result:

\begin{lemma}
\label{lemma:half-form density limit}Let $f\in C^{\infty}(M)$ and define%
\[
J_{k,R}(f)(x_{0})=\left(  k/2\pi\right)  ^{d/2}\int_{B_{R}(0)}e^{-k\rho
(\xi,x_{0})}f(\xi,x_{0})\exp\left\{  -\int_{\gamma_{\xi}}\frac{{\mathcal{L}%
}_{JX^{\xi}}\varepsilon_{\omega}}{2\varepsilon_{\omega}}\right\}  \tau
(\xi,x_{0})d^{d}\xi.
\]
Then%
\[
\lim_{k\rightarrow\infty}\left\vert J_{k,R}(x_{0})-f(x_{0})\right\vert =0
\]
uniformly on $\Phi^{-1}(0).$
\end{lemma}

We are now ready to prove our main result.

\begin{proof}[Proof of Theorem \protect\ref{thm:density asymptotics}]
We write $I_{k}$ as the sum of an integral over $B_{R}(0)$ and an integral
over the complement of $B_{R}(0)$. Combining Theorem \ref{thm:s along a curve} and
Lemmas \ref{lemma:tau=1} and \ref{lemma:density limit}, the first term
approaches $2^{-d/2}f(x_{0})\mathrm{\operatorname{vol}}(G\cdot x_{0})$ uniformly
as $k\rightarrow \infty $. By Lemma \ref{lemma:negligable contributions} the
second term approaches $0$ uniformly as $k\rightarrow \infty $.
If we include the metaplectic correction, the proof is similar, except the
first term approaches $f(x_{0})$ uniformly.
\end{proof}

\section{Discussion and examples}

We begin by noting that compactness is not essential to the issues considered
in this paper. Given a noncompact K\"{a}hler manifold $M$ acted on by a
compact group $G,$ one can define a natural map between the $G$-invariant
holomorphic sections of the relevant line bundle over $M$ and the space of all
holomorphic sections of the ``quotient'' bundle over $M/\!\!/G.$ (Simply restrict the section over $M$ to
$\Phi^{-1}(0)$ and then let it descend to $M/\!\!/G=\Phi^{-1}(0)/G.$) Although
it is unlikely that this map is invertible for arbitrary noncompact $M,$ it is
likely to be invertible in many interesting examples, and similarly in the
presence of the metaplectic correction. It therefore makes sense to
investigate the issue of unitarity at least in the more favorable noncompact examples.

Indeed, it even makes sense, to some extent, to consider quantization and
reduction for certain \textit{infinite-dimensional} K\"{a}hler manifolds. This
sort of problem arises naturally in quantum field theories, where the
holomorphic approach to quantization is often the most natural one and where
one usually has to reduce by an (infinite-dimensional) group of gauge
symmetries. Of course, it inevitably requires some creativity to give a
sensible meaning to the quantization and to the reduction in
infinite-dimensional settings. Nevertheless, there are some interesting
examples (discussed below) where this can be done.

In the rest of this subsection, we discuss some noncompact (and, in some
cases, infinite dimensional) examples in which the issue of unitarity in
quantization versus reduction is of interest. In some of these cases, a
Guillemin--Sternberg-type map (with the metaplectic correction) actually turns
out to be exactly unitary.

Our first example is the quantization of Chern--Simons theory, as considered
in the paper \cite{Axelrod} of Axelrod, della Pietra, and Witten. The authors
perform a K\"{a}hler quantization of the moduli space $\mathcal{M}$ of flat
connections modulo gauge transformations over a Riemann surface
$\mathrm{\Sigma}.$ Much of the analysis is done by regarding this moduli space
as the symplectic quotient of the space $\mathcal{A}$ of all connections by
the group $\mathcal{G}$ of gauge transformations, using the result of Atiyah
and Bott that the moment map for the action of $\mathcal{G}$ is simply the
curvature. The main result of \cite{Axelrod} is the construction of a natural
projectively flat connection that serves to identify the Hilbert spaces
obtained by using different complex structures on $\mathrm{\Sigma}.$ The
existence of the projectively flat connection shows that the quantization
procedure is independent of the choice of complex structure on $\mathrm{\Sigma
},$ since it allows one to identify (projectively) all the different Hilbert
spaces with one another.

There is, however, one important issue that is not fully resolved in
\cite{Axelrod}, namely the issue of the unitarity of the connection. The
projectively flat connection \textquotedblleft upstairs\textquotedblright\ on
$\mathcal{A}$ is unitary, at least formally. The authors suggest, then, that
one should simply define the norm of a section downstairs to be the norm of
the corresponding section upstairs, in which case, the connection would
(formally) be unitary. However, because $\mathcal{A}$ is infinite-dimensional,
it is not entirely clear that this prescription makes sense. Now, \textit{if}
it were true that the Guillemin--Sternberg map was unitary in this context,
that would mean that computing the norm upstairs is the same as computing the
norm downstairs. In that case, one would expect to have unitarity using the
natural downstairs norm on the space of sections. Since (as we show in this
paper) one cannot expect the Guillimen--Sternberg map to be unitary in
general, it is not clear what norm one should use in order to have the
projectively flat connection be unitary. Thus, the failure of unitarity (in
general) has consequences in this case.

A second example is the canonical quantization of $(1+1)$-dimensional
Yang--Mills theory. In this case, the upstairs space is an
infinite-dimensional affine space, namely the cotangent bundle of the space of
connections over the spatial circle. Because the upstairs space is an affine
space, there is a well-defined (K\"{a}hler) quantization, namely a
Segal--Bargmann space over an infinite-dimensional vector space (as in
\cite{Baez}). Unfortunately, as is often the case in field theories, there are
no nonzero vectors in this space that are invariant under the action of the
gauge group \cite{DHeRep}. Thus, if one wants to do quantization first and
then reduction, some regularization procedure must be used when performing the
reduction. Two different regularization procedures have been considered, that
of Wren \cite{LandsmanWren},\cite{Wren} (using Landsman's generalized Rieffel
induction \cite{Landsman}) and that of Driver--Hall \cite{DHym} (using a
Gaussian measure of large variance to approximate the nonexistent Lebesgue
measure). The two procedures give the same result, that the \textquotedblleft
first quantize and then reduce\textquotedblright\ space can be identified with
a certain $L^{2}$-space of holomorphic functions over the complexified
structure group $G_{\mathbb{C}}.$

Meanwhile, the results of \cite{HallGQ} indicate that the same $L^{2}$-space
of holomorphic functions can be obtained by geometric quantization of
$G_{\mathbb{C}},$ provided that the metaplectic correction is used. This means
that in this case, a Guillemin--Sternberg-type map, with metaplectic
correction, does turn out to be unitary. See also \cite{HallCSYM}.

A third example, related to the second one, is reduction from $G_{\mathbb{C}}
$ (which is naturally identified with the cotangent bundle $T^{\ast}(G)$) to
$G_{\mathbb{C}}/H_{\mathbb{C}}$ (identified with $T^{\ast}(G/H)$). Here $H $
is the fixed-point subgroup of an involution, which means that $G/H$ is a
compact symmetric space. Results of Hall \cite{HallGQ} and Stenzel
\cite{Stenzel} show that the \textquotedblleft first quantize and then
reduce\textquotedblright\ space can be identified with an $L^{2}$-space of
holomorphic functions on $G_{\mathbb{C}}/H_{\mathbb{C}}$ with respect to a
certain heat kernel measure. (This result can also be obtained from the
computation of the relevant orbital integral by Flensted-Jensen \cite[Eq.
6.20]{FJ}.) \textit{If} $G/H$ is again isometric to a compact Lie group (i.e.,
if $G=U\times U$ and $H$ is the diagonal copy of $U$), then the
above-mentioned results show that the \textquotedblleft first reduce and then
quantize\textquotedblright\ can be identified in a natural unitary fashion
with the \textquotedblleft first quantize and then reduce\textquotedblright%
\ space, provided that the metaplectic correction is included in both
constructions. So this situation provides another example in which a
Guillemin--Sternberg-type map (with metaplectic correction) is unitary.

On the other hand, if $G/H$ is not isometric to a compact Lie group, then the
norms on the two spaces are not the same. Rather, the measure obtained in the
\textquotedblleft first reduce and then quantize\textquotedblright\ space is
the leading term in the asymptotic expansion of the heat kernel measure in the
\textquotedblleft first quantize and then reduce\textquotedblright\ space; see
the discussion on pp. 244--245 of \cite{HallGQ}. (In the case that $G/H$ is
isometric to a compact Lie group, the relevant heat kernel is actually equal
to the leading term in the expansion, up to an overall constant.) In this
case, we see that, even with the metaplectic correction, the measures used to
compute the two norms are not equal; they are, however, equal to leading order
in $\hbar.$

A final example, considered in \cite{FMMN}, is reduction from $T^{\ast}(G)$ to
$T^{\ast}(G/\mathrm{\operatorname*{Ad}}G).$ We may identify
$G/\mathrm{\operatorname*{Ad}}G$ with $T/W,$ where $T$ is a maximal torus in
$G$ and $W$ is the Weyl group. Although in this case the action of $G$ is
generically nonfree, one should be able to construct a natural map (with the
metaplectic correction) by a construction similar to the one we give in this
paper. We expect that this map will turn out to be the one given by
$F\rightarrow\sigma\left.  F\right\vert _{T}$ in the notation of \cite{FMMN}.
If this is the case, then Theorems 2.2 and 2.3 of \cite{FMMN} will show that
the identification of the \textquotedblleft first quantize and then reduce
space\textquotedblright\ with the \textquotedblleft first reduce and then
quantize\textquotedblright\ space is unitary, up to a constant. The results of
\cite{FMMN} are related to the genus-one case of the quantization of
Chern--Simons theory considered in \cite{Axelrod}.

\subsection*{Acknowledgements}

The authors are grateful to Dennis Snow for suggesting to use an embedding
into projective space to prove Lemma \ref{lemma:tau estimates}. They would
also like to thank Alejandro Uribe for bringing to their attention the papers
\cite{Charles} and \cite{Charles2}, and Liviu Nicolaescu for pointing us
toward references for the coarea formula.

\providecommand{\bysame}{\leavevmode\hbox to3em{\hrulefill}\thinspace}

\end{document}